\documentclass{gtart}
\def\ifplaintex{\expandafter\ifx\csname documentclass\endcsname\relax}

\ifplaintex 
\else
\fi

\expandafter\ifx\csname beginpicture\endcsname\relax
\expandafter\ifx\csname documentclass\endcsname\relax
\input pictex \else
\input prepictex \input pictex \input postpictex \fi\fi

\def\gt{{\mathsurround=0pt\it $\cal G\mskip-2mu$eometry \&\ 
$\cal T\!\!$opology}}        %  journal title in recommended style

\def\gtp{{\mathsurround=0pt\it $\cal G\mskip-2mu$eometry \&\ 
$\cal T\!\!$opology $\cal P\!$ublications}}  % GT publications

\def\lognumber#1{\def\thelognumber{#1}}
\def\volumenumber#1{\def\thevolumenumber{#1}}
\def\papernumber#1{\def\thepapernumber{#1}}
\def\volumeyear#1{\def\thevolumeyear{#1}}

\def\pagenumbers#1#2{\def\startpage{#1}\def\finishpage{#2}}
\def\published#1{\def\publishdate{#1}}
\def\proposed#1{\def\theproposer{#1}}
\def\seconded#1{\def\theseconders{#1}}
\def\received#1{\def\receiveddate{#1}}

\def\accepted#1{\def\accepteddate{#1}}

\def\coverauthors#1{\def\thecoverauthors{#1}}
\def\asciiauthors#1{\def\theasciiauthors{#1}}

\long\def\asciiabstract#1{\long\def\theasciiabstract{#1}}
\def\asciikeywords#1{\def\theasciikeywords{#1}}

\def\shorttitle#1{\def\theshorttitle{#1}}

\let\\\par\let\thelognumber\relax
\let\thevolumenumber\relax\let\thepapernumber\relax
\let\thevolumeyear\relax\let\thesamplenumber\relax\let\startpage\relax
\let\finishpage\relax\let\publishdate\relax\let\receiveddate\relax
\let\reviseddate\relax\let\accepteddate\relax\let\theasciititle\relax
\let\theasciiauthors\relax
\let\theasciiabstract\relax\let\theasciikeywords\relax
\let\theasciiemail\relax\let\theshortauthors\relax\let\theshorttitle\relax
\let\thecoverauthors\relax

\long\def\maketitlep{   % start of definition of \maketitlep

\count0=\startpage

\gt\hfill      %   Journal title (top left) 
\beginpicture
\setcoordinatesystem units <0.33truein, 0.33truein> point at 2.2 0.9
\setplotsymbol ({$\cal G$})
\plotsymbolspacing=9truept
\circulararc 315 degrees from 0 1 center at 0 0
\setplotsymbol ({$\cal T$})
\circulararc 315 degrees from 1 -1 center at 1 0
\endpicture
\break
{\small\ifx\thesamplenumber\relax % sample?  
Volume \else Sample
\fi\thevolumenumber\ (\thevolumeyear)
\startpage--\finishpage\nl
Published: \publishdate}
\vglue 0.5truein plus 0.4fil minus 0.1truein

{\parskip=0pt\leftskip 0pt plus 1fil\def\\{\par\smallskip}{\ifplaintex\large
\else\Large\fi\bf\thetitle}\par\medskip}   

\vglue 0pt plus 0.1fil 

{\parskip=0pt\leftskip 0pt plus 1fil\def\\{\par}{\sc\theauthors}
\par\medskip}

\vglue 0pt plus 0.1fil 

{\small\parskip=0pt\let\newline\\
{\leftskip 0pt plus 1fil\def\\{\par}{\sl\theaddress}\par}
\expandafter\ifx\theemail\relax    % email address?
\relax\else\vglue 5pt plus 0.02fil minus 2pt\def\\{\stdspace{\rm 
and}\stdspace} 
\cl{Email:\stdspace\tt\theemail}\fi
\ifx\theurl\relax                  % URL given?
\relax\else\vglue 5pt plus 0.02fil minus 2pt\def\\{\stdspace{\rm 
and}\stdspace}
\cl{URL:\stdspace\tt\theurl}\fi\par}

\vglue 7pt plus 0.3fil minus 3pt

{\bf Abstract}
\vglue 5pt plus 0.1fil minus 2pt

\theabstract

\vglue 7pt plus 0.3fil minus 3pt

{\bf AMS Classification numbers}\quad Primary:\quad \theprimaryclass

Secondary:\quad \thesecondaryclass

\vglue 5pt plus 0.3fil minus 2pt

{\bf Keywords}\quad \thekeywords

\vglue 10pt plus 0.5fil minus 5pt

{\small  Proposed: \theproposer\hfill Received: \receiveddate\nl
Seconded: \theseconders\hfill 
\ifx\reviseddate\relax                         % paper revised?
Accepted: \accepteddate                        % no
\else
Revised: \reviseddate                          % yes
\fi}
\eject
}       %  end of definition of \maketitlep

\let\maketitlepage\maketitlep
\let\maketitle\maketitlepage

\font\phead=cmsl9 scaled 950
\font=cmsl9 scaled 1050
\font\pnum=cmbx10 scaled 913
\font\lnum=cmbx10 
\font\pfoot=cmsl9 scaled 950
\font=cmsl9 scaled 1050
\ifplaintex
\headline{\vbox to 0pt{\vskip -4.5mm\line{\small\phead\ifnum
\count0=\startpage ISSN 1364-0380 (on line)
1465-3060 (printed) \hfill {\pnum\folio}\else\ifodd\count0\def\\{ }% 
\ifx\theshorttitle\relax\thetitle\else\theshorttitle\fi\hfill{\pnum\folio}
\else\def\\{ and }{\pnum\folio}\hfill\ifx\theshortauthors\relax\theauthors
\else\theshortauthors\fi\fi\fi}\vss}}
\footline{\vbox to 0pt{\vglue 0mm\line{\small\pfoot\ifnum\count0=\startpage
\copyright\ \gtp\hfill\else
\gt, Volume \thevolumenumber\ (\thevolumeyear)\hfill\fi}\vss
}}
\else
\makeatletter
\def\@oddhead{{\small\lhead\ifnum\count0=\startpage ISSN 1364-0380 (on line)
1465-3060 (printed) \hfill {\lnum\number\count0}\else\ifodd\count0
\def\\{ }\ifx\theshorttitle\relax \thetitle \else\theshorttitle\fi\hfill
{\lnum\number\count0}\else\def\\{ and }{\lnum\number\count0}
\hfill\ifx\theshortauthors\relax 
\theauthors\else\theshortauthors\fi\fi\fi}}\def\@evenhead{\@oddhead}
\def\@oddfoot{\small\lfoot\ifnum\count0=\startpage\copyright\ \gtp\hfill\else
\gt, Volume \thevolumenumber\ (\thevolumeyear)\hfill\fi}
\def\@evenfoot{\@oddfoot}
\makeatother
\fi

\newwrite\gtoutfile
\long\gdef\makeheadfile{  %%% start of definition of \makeheadfile
{\def\\{, }\def\s{ }
\immediate\openout\gtoutfile head.xxx
\immediate\write\gtoutfile{Proxy-for: \ifx\theasciiauthors\relax
\theauthors\else\theasciiauthors\fi\s<\ifx\theasciiemail\relax\theemail\else\theasciiemail\fi>}
\immediate\write\gtoutfile{\noexpand\\}
\immediate\write\gtoutfile{Authors: \ifx\theasciiauthors\relax
\theauthors\else\theasciiauthors\fi}
{\def\\{ }\immediate\write\gtoutfile{Title: \ifx\theasciititle\relax
\thetitle\else\theasciititle\fi}}
\immediate\write\gtoutfile{Subj-class: GT or SG or MG etc}
\immediate\write\gtoutfile{MSC-class: \theprimaryclass\ifx\thesecondaryclass\relax\else, \thesecondaryclass\fi}
\immediate\write\gtoutfile{Journal-ref: Geom. Topol. \thevolumenumber
(\thevolumeyear) \startpage-\finishpage}
\immediate\write\gtoutfile{Comments: Published by Geometry and Topology at}
\immediate\write\gtoutfile{\s\s http://www.maths.warwick.ac.uk/gt/GTVol\thevolumenumber/paper\thepapernumber.abs.html}
\immediate\write\gtoutfile{\noexpand\\}
\immediate\write\gtoutfile{}
\ifx\theasciiabstract\relax
\immediate\write\gtoutfile{\theabstract}\else
\immediate\write\gtoutfile{\theasciiabstract}\fi
\immediate\write\gtoutfile{}
\immediate\write\gtoutfile{\noexpand\\}
\immediate\write\gtoutfile{}
\immediate\closeout\gtoutfile}}  %%% end of definition of \makeheadfile

\def\maketitlepage{\maketitlep\makeheadfile}
\let\maketitle\maketitlepage

\lognumber{386}
\volumenumber{7}\papernumber{29}\volumeyear{2003}
\pagenumbers{1001}{1053}
\received{26 November 2003}
\published{26 December 2003}
\accepted{25 December 2003}
\proposed{Haynes Miller}
\seconded{Thomas Goodwillie, Gunnar Carlsson}

\usepackage{amscd, amsmath}
\usepackage{amssymb}

\newcommand{\tensor}{\otimes}
\newcommand{\homeq}{\simeq}
\newcommand{\iso}{\cong}

\DeclareMathOperator{\Ho}{Ho}
\DeclareMathOperator*{\colim}{colim}
\DeclareMathOperator{\map}{map}
\DeclareMathOperator{\Hom}{Hom}
\DeclareMathOperator{\Tot}{Tot}
\DeclareMathOperator{\sk}{sk}
\DeclareMathOperator{\Fib}{Fib}
\DeclareMathOperator{\diag}{diag}
\DeclareMathOperator{\con}{con}
\DeclareMathOperator{\core}{core}
\DeclareMathOperator{\Cotor}{Cotor}
\DeclareMathOperator{\Ext}{Ext}

\newcommand{\Z}{\mathbb{Z}}
\newcommand{\B}{\mathcal{B}}
\newcommand{\C}{\mathcal{C}}
\newcommand{\D}{\mathcal{D}}
\newcommand{\G}{\mathcal{G}}
\newcommand{\Hcal}{\mathcal{H}}
\newcommand{\I}{\mathcal{I}}
\newcommand{\Lcal}{\mathcal{L}}
\newcommand{\M}{\mathcal{M}}
\newcommand{\N}{\mathcal{N}}
\newcommand{\R}{\mathcal{R}}
\newcommand{\Scal}{\mathcal{S}}
\newcommand{\cC}{c\,\mathcal{C}}
\newcommand{\cD}{c\,\mathcal{D}}
\newcommand{\cM}{c\,\mathcal{M}}
\newcommand{\cS}{c\,\mathcal{S}}

\newcommand{\Deltabf}{\boldsymbol{\Delta}}

\makeatletter
\let\c@equation\c@subsection

\makeatother

\theoremstyle{plain} 
\newtheorem{thm}[equation]{Theorem}
\newtheorem{cor}[equation]{Corollary}
\newtheorem{lem}[equation]{Lemma}
\newtheorem{prop}[equation]{Proposition}

\theoremstyle{definition}
\newtheorem{defn}[equation]{Definition}

\newtheorem{ex}[equation]{Example}

\begin{document}

\title{Cosimplicial resolutions and homotopy spectral\\sequences in 
model categories}
\shorttitle{Cosimplicial resolutions and homotopy spectral sequences}
\author{A\,K Bousfield}
\coverauthors{A\noexpand\thinspace K Bousfield}
\asciiauthors{A K Bousfield}
\address{Department of Mathematics\\
         University of Illinois at Chicago\\
         Chicago, Illinois 60607, USA}

\email{bous@uic.edu}

\asciikeywords{Cosimplicial resolutions, homotopy spectral sequences, 
model categories, Bendersky-Thompson completion, Bousfield-Kan completion}

\begin{abstract}
We develop a general theory of cosimplicial resolutions, homotopy
spectral sequences, and completions for objects in model categories,
extending work of Bousfield--Kan and Bendersky--Thompson for ordinary
spaces.  This is based on a generalized cosimplicial version of the
Dwyer--Kan--Stover theory of resolution model categories, and we are
able to construct our homotopy spectral sequences and completions
using very flexible weak resolutions in the spirit of relative
homological algebra.  We deduce that our completion functors have
triple structures and preserve certain fiber squares up to homotopy.
We also deduce that the Bendersky--Thompson completions over connective
ring spectra are equivalent to Bousfield--Kan completions over solid
rings.  The present work allows us to show, in a subsequent paper,
that the homotopy spectral sequences over arbitrary ring spectra have
well-behaved composition pairings.
\end{abstract}
\asciiabstract{%
We develop a general theory of cosimplicial resolutions, homotopy
spectral sequences, and completions for objects in model categories,
extending work of Bousfield-Kan and Bendersky-Thompson for ordinary
spaces.  This is based on a generalized cosimplicial version of the
Dwyer-Kan-Stover theory of resolution model categories, and we are
able to construct our homotopy spectral sequences and completions
using very flexible weak resolutions in the spirit of relative
homological algebra.  We deduce that our completion functors have
triple structures and preserve certain fiber squares up to homotopy.
We also deduce that the Bendersky-Thompson completions over connective
ring spectra are equivalent to Bousfield-Kan completions over solid
rings.  The present work allows us to show, in a subsequent paper,
that the homotopy spectral sequences over arbitrary ring spectra have
well-behaved composition pairings.}

\maketitlepage

\section{Introduction}
\label{sec1}

In \cite{BK} and \cite{BK1972a}, Bousfield--Kan developed unstable Adams spectral sequences and completions of spaces with respect to a ring, and this work was extended by Bendersky--Curtis--Miller \cite{BCM} and Bendersky--Thompson \cite{BT} to allow a ring spectrum in place of a ring. In the present work, we develop a much more general theory of cosimplicial resolutions, homotopy spectral sequences, and completions for objects in model categories.  Among other things, this provides a flexible approach to the Bendersky--Thompson spectral sequences and completions, which is especially needed because the original chain level constructions of pairings and products \cite{BK1973} do not readily extend to that setting.

We rely heavily on a generalized cosimplicial version of the Dwyer--Kan--Stover \cite{DKS} theory of \emph{resolution model categories} (or \emph{$E_2$ model categories} in their parlance).  This provides a simplicial model category structure $\cC^{\G}$ on the category $\cC$ of cosimplicial objects over a left proper model category $\C$ with respect to a chosen class $\G$ of \emph{injective models} (see Theorems \ref{thm:3.3} and \ref{thm:12.4}).  Of course, our cosimplicial statements have immediate simplicial duals.  Other more specialized versions of the simplicial theory are developed by Goerss--Hopkins \cite{GH} and Jardine \cite{Jar} using small object arguments which are not applicable in the duals of many familiar model categories. When $\C$ is discrete, our version reduces to a variant of Quillen's model category structure \cite[II\S4]{Qui} on $c\C$, allowing many possible choices of ``relative injectives'' in addition to Quillen's canonical choice (see \ref{sec:4.3} and \ref{sec:4.4}).  However, we are most interested in examples where $\C$ is the category of pointed spaces and where $\G$ is determined by a ring spectrum (\ref{sec:4.9}) or a cohomology theory (\ref{sec:4.6}).  In the former case, the model category provides Bendersky--Thompson-like \cite{BT} cosimplicial resolutions of spaces with respect to an arbitrary ring spectrum, which need not be an $S$--algebra.

In general, a \emph{cosimplicial $\G$--injective resolution}, or \emph{$\G$--resolution}, of an object $A\in\C$ consists of a trivial cofibration $A\to\bar{A}^\bullet$  to a fibrant target $\bar{A}^\bullet$ in $\cC^{\G}$.  By applying the constructions of \cite{BK} and \cite{DP} to $\G$--resolutions, we obtain \emph{right derived functors} $\R^s_{\G}T(A)=\pi^sT(\bar{A}^\bullet)$, \emph{$\G$--completions} $\hat{L}_{\G}A=\Tot\bar{A}^\bullet$, and \emph{$\G$--homotopy spectral sequences} 
$\{E^{s,t}_r(A;M)_{\G}\}_{r\geq 2}=\{E^{s,t}_r(\bar{A}^\bullet;M)\}_{r\geq 2}$  abutting to $[M,\hat{L}_{\G}A]_*$ for $A,M\in\C$  (see \ref{sec:5.5}, \ref{sec:5.7}, and \ref{sec:5.8}).  We proceed to show that the $\G$--resolutions in these constructions may be replaced by \emph{weak $\G$--resolutions}, that is, by arbitrary weak equivalences in $\cC^{\G}$ to termwise $\G$--injective targets (see Theorems \ref{thm:6.2} and \ref{thm:6.5}).  This is convenient since weak $\G$--resolutions are easy to recognize and arise naturally from triples on $\C$.  The Bendersky--Thompson resolutions are clearly examples of them.

We deduce that the $\G$--completion functor $\hat{L}_{\G}$ belongs to a triple on the homotopy category $\Ho\C$ (see Corollary \ref{cor:8.2}), and we introduce notions of \emph{$\G$--completeness}, \emph{$\G$--goodness}, and \emph{$\G$--badness} for objects in $\Ho\C$.  This generalizes work of Bousfield--Kan \cite{BK} on the homotopical $R$--completion functor $R_\infty$ for pointed spaces.  We discuss an apparent error in the space-level associativity part of the original triple lemma \cite[page 26]{BK} for $R_\infty$, but we note that this error does not seem to invalidate any of our other results (see \ref{sec:8.9}).  We also develop criteria for comparing different completion functors, and we deduce that the Bendersky--Thompson completions with respect to connective ring spectra are equivalent to Bousfield--Kan completions with respect to solid rings (see Theorem \ref{thm:9.8}), even though the associated homotopy spectral sequences may be very different.

Finally, we show that the $\G$--completion functors preserve certain fiber squares up to homotopy (see Theorem \ref{thm:10.9}), and we focus particularly on the Bendersky--Thompson $K$--completion and the closely related $p$--adic $K$--completion, where $K$ is the spectrum of nonconnective $K$--theory at a prime $p$.  In particular, we find that the $K$--completion functor preserves homotopy fiber squares when their $K_*$--cobar spectral sequences collapse strongly and their spaces have free $K_*$--homologies, while the $p$--adic $K$--completion functor preserves homotopy fiber squares when their $K/p_*$--cobar spectral sequences collapse strongly and their spaces have torsion-free $p$--adic $K$--cohomologies (see Theorems \ref{thm:10.12} and \ref{thm:11.7}).  In general, the $K$--completions and $K$--homotopy spectral sequences are very closely related to their $p$--adic variants (see Theorem \ref{thm:11.4}), though the latter seem to have better technical properties.  For instance, the $p$--adic $K$--homotopy spectral sequences seem especially applicable to spaces whose $p$--adic $K$--cohomologies are torsion-free with Steenrod--Epstein-like $U(M)$ structures as in \cite{Bou1996a}.

In much of this work, for simplicity, we assume that our model categories are pointed.  However, as in \cite{GH}, this assumption can usually be eliminated, and we offer a brief account of the unpointed theory in Section \ref{sec:12}.  We thank Paul Goerss for suggesting such a generalization.

In a sequel \cite{Bou2002}, we develop composition pairings for our homotopy spectral sequences and discuss the $E_2$--terms from the standpoint of homological algebra.  This extends the work of \cite{BK1973}, replacing the original chain-level formulae over rings by more general constructions.  It applies to give composition pairings for the Bendersky--Thompson spectral sequences.

Although we have long been interested in the present topics, we were prompted to formulate this theory by Martin Bendersky and Don Davis who are using some of our results in \cite{BD} and \cite{BD2001}, and we thank them for their questions and comments.  We also thank Assaf Libman for his suggestions and thank the organizers of BCAT 2002 for the opportunity to present this work. 

Throughout, we assume a basic familiarity with Quillen model categories and generally follow the terminology of \cite{BK}, so that ``space'' means ``simplicial set.''  The reader seeking a rapid path into this work might now review the basic terminology in Section \ref{sec:2}, then read the beginning of Section \ref{sec:3} through the existence theorem (\ref{thm:3.3}) for resolution model categories, and then proceed to the discussion of these categories in Section \ref{sec:4}, skipping the very long existence proof in Section \ref{sec:3}.

 The paper is divided into the following sections:
\begin{enumerate}
\item[1.]Introduction
\item[2.]Homotopy spectral sequences of cosimplicial objects
\item[3.]Existence of resolution model categories
\item[4.]Examples of resolution model categories
\item[5.]Derived functors, completions, and homotopy spectral sequences
\item[6.]Weak resolutions are sufficient
\item[7.]Triples give weak resolutions
\item[8.]Triple structures of completions
\item[9.]Comparing different completions
\item[10.]Bendersky--Thompson completions of fiber squares
\item[11.]$p$--adic $K$--completions of fiber squares
\item[12.]The unpointed theory
\end{enumerate}

The author was partially supported by the National Science Foundation.

\section{Homotopy spectral sequences of cosimplicial objects}
\label{sec:2}

We now introduce the homotopy spectral sequences of cosimplicial objects in model categories, thereby generalizing the constructions of Bousfield--Kan \cite{BK} for cosimplicial spaces.  This generalization is mainly due to Reedy\cite{Ree}, but we offer some details to establish notation and terminology.  We first consider the following:

\subsection{Model categories}
\label{sec:2.1}
By a \emph{model category} we mean a closed model category in Quillen's original sense \cite{Qui}.  This consists of a category with three classes of maps called \emph{weak equivalences}, \emph{cofibrations}, and \emph{fibrations}, satisfying the usual axioms labeled {\bf MC1}--{\bf MC5} in \cite[pages 83--84]{DS}.  We refer the reader to \cite{DS}, \cite{GJ}, \cite{Hir}, and \cite{Hov} for good recent treatments of model categories.  A model category is called \emph{bicomplete} when it is closed under all small limits and colimits.  It is called \emph{factored} when the factorizations provided by {\bf MC5} are functorial.  We note that most interesting model categories are bicomplete and factored or factorable, and some authors incorporate these conditions into the axioms (see \cite{Hir} and \cite{Hov}).

\subsection{Cosimplicial objects}
\label{sec:2.2}
A \emph{cosimplicial object} $X^\bullet$ over a category $\C$ consists of a diagram in $\C$ indexed by the category    $\Deltabf$ of finite ordinal numbers.  More concretely, it consists of objects $X^n \in \C$ for $n \geq 0$ with \emph{coface} maps $d^i\co X^n \to X^{n+1}$ for $0 \leq i \leq{n+1}$ and \emph{codegeneracy} maps $s^i\co X^{n+1}\to X^n$ for $0\leq j \leq n$ satisfying the usual cosimplicial identities (see \cite[page 267]{BK}).  Thus a cosimplicial object over $\C$ corresponds to a simplicial object over $\C ^{\rm op}$.  The category of cosimplicial objects over $\C$ is denoted by $\cC$, while that of simplicial objects is denoted by $s\,\C$.  

When $\C$ is a model category, there is an induced model category structure on $\cC =s(\C ^{\rm op})$ due to Reedy~\cite{Ree}.  This is described by Dwyer--Kan--Stover~\cite{DKS}, Goerss--Jardine~\cite{GJ}, Hirschhorn~\cite{Hir}, Hovey~\cite{Hov}, and others.  For an object $X^\bullet\in \cC$, consider the \emph{latching} maps $L^n X^\bullet\to X^n$ in $\C$ for $n\geq 0$ where
$$L^n X^\bullet = \colim_{\theta\co [k] \to [n]}X^k$$ 
with $\theta$ ranging over the injections  $[k]\to [n]$  in $\Deltabf$ for 
$k<n$, and consider the \emph{matching} maps $X^n\to M^n X^\bullet$ in $\C$ for   $n\geq 0$ where 
$$M^n X^\bullet = \lim_{\phi \co [n]\to [k]}X^k$$
with $\phi$ ranging over the surjections $[n]\to [k]$ in $\Deltabf$ for $k<n$.  A cosimplicial map $f\co X^\bullet\to Y^\bullet \in \cC$ is called:
\begin{enumerate}
\item[(i)]a \emph{Reedy weak equivalence} when $f\co X^n\to Y^n$ is a weak equivalence in $\C$ for $n\geq 0$;
\item[(ii)]a \emph{Reedy  cofibration} when $X^n\coprod_{L^n X^\bullet}L^n Y^\bullet\to Y^n$ is a cofibration in $\C$ for $n\geq 0$;
\item[(iii)]a \emph{Reedy fibration} when $X^n\to Y^n\times_{M^n Y^\bullet}M^n X^\bullet$ is a fibration in $\C$ for $n\geq 0$.
\end{enumerate}

\begin{thm}[Reedy]
\label{thm:2.3}
If $\C$ is a model category, then so is $\cC$ with the Reedy weak equivalences, Reedy cofibrations, and Reedy fibrations.
\end{thm}

\begin{ex}
\label{ex:2.4}
Let $\Scal$ and $\Scal_*$ denote the categories of spaces (ie, simplicial sets) and pointed spaces with the usual model category structures.  Then the Reedy model category  structures on $\cS$ and $\cS_*$ reduce to those of Bousfield--Kan \cite[page 273]{BK}.  Thus a map $X^\bullet\to Y^\bullet$ in $\cS$ or $\cS_*$ is a Reedy weak equivalence when it is a termwise weak equivalence, and is a Reedy cofibration when it is a termwise injection such that $a(X^\bullet)\iso a(Y^\bullet)$ where $a(X^\bullet)=\{ x\in X^0~|~d^0x=d^1x\}$ is the maximal augmentation.
\end{ex}

\subsection{Simplicial model categories}
\label{sec:2.5}
As in Quillen \cite[II.1]{Qui}, by a \emph{simplicial category},  we mean a category $\C$ enriched over $\Scal$, and we write $\map(X,Y) \in \Scal$ for the mapping space of $X,Y\in\C$.  When they exist, we also write $X\tensor K\in\C$ and $\hom (K,X)\in\C$ for the tensor and cotensor of $X\in\C$ with $K\in\Scal$.  Since there are natural equivalences
$$\Hom_{\Scal}(K,\map(X,Y))~\iso~\Hom_{\C}(X\tensor K,Y)~\iso~\Hom_{\C}(X,\hom(K,Y)),$$
any one of the three functors, $\map$, $\tensor$, and $\hom$, determines the other two uniquely.  As in Quillen \cite[II.2]{Qui}, by a \emph{simplicial model category},  we mean a model category $\C$ which is also a simplicial category satisfying the following axioms {\bf SM0} and {\bf SM7} (or equivalently {\bf SM7$^\prime$}):
\begin{description}
\item[SM0]The objects $X\tensor K$ and $\hom(K,X)$ exist for each $X\in\C$ and each finite $K\in\Scal$.
\item[SM7] If $i\co A\to B\in\C$ is a cofibration and $p\co X\to Y\in\C$ is a fibration, then the map
$$\map(B,X)\xrightarrow{\qquad}\map(A,X)\times_{\map(A,Y)}\map(B,Y)$$
is a fibration in $\Scal$ which is trivial if either $i$ or $p$ is trivial. 
\item[SM7$^\prime$]If $i\co A\to B\in\C$ and $j\co J\to K\in\Scal$ are cofibrations with $J$ and $K$ finite, then the map
$$(A\tensor K)\coprod_{A\tensor J}(B\tensor J)\xrightarrow{\qquad} B\tensor K$$
is a cofibraton in $\C$ which is trivial if either $i$ or $j$ is trivial.
\end{description}

\begin{thm}
\label{thm:2.6}
If C is a simplicial model category, then so is the Reedy model category $\cC$ with $(X^\bullet\tensor K)^n = X^n\tensor K$ and $\hom(K,X^\bullet)^n=\hom(K,X^n)$ for $X^\bullet\in\cC$ and finite $K\in\Scal$.
\end{thm}

\begin{proof}
The simplicial axiom {\bf SM7$^\prime$} follows easily using the isomorphisms\break $L^n(X^\bullet\tensor K)\iso L^n X^\bullet\tensor K$ for $n\geq 0$.
\end{proof}

To construct our total objects and spectral sequences, we need the following:

\subsection{Prolongations of the mapping functors}
\label{sec:2.7}
Let $\C$ be a bicomplete simplicial model category.  Then the objects $X\tensor K\in\C$ and $\hom(K,X)\in\C$ exist for each $X\in\C$ and each $K\in \Scal$, without finiteness restrictions.  For $A\in\C$, $Y^\bullet\in\cC$, and $J^\bullet\in\cS$, we define $\map(A,Y^\bullet)\in\cS$ and $A\tensor J^\bullet\in\cC$ termwise, and we let \ $\hom(J^\bullet,-)\co \cC\to\C$ \ denote the right adjoint of \ 
$-\tensor J^\bullet\co \C\to\cC$.  It is not hard to show that the functor \ $\tensor\co \C\times\cS\to\cC$ \ satisfies the analogue of {\bf SM7$^\prime$}, and hence the functors  \ $\map\co \C^{\rm op}\times\cC\to\cS$ \  and \  
$\hom\co (\cS)^{\rm op}\times\cC\to \C$ \  satisfy the analogues of {\bf SM7}.

\subsection{Total objects}
\label{sec:2.8}
Now let $\C$ be a pointed bicomplete simplicial model category, and let $X^\bullet\in\cC$ be Reedy fibrant.  The \emph{total object}  $\Tot X^\bullet = \hom(\Delta^\bullet,X^\bullet)\in\C$ is defined using the cosimplicial space $\Delta^\bullet\in\cS$ of standard $n$--simplices $\Delta^n\in\Scal$ for $n\geq 0$.  It is the limit of the \emph{Tot tower} $\{ \Tot_s X^\bullet\} _{s\geq 0}$ with $\Tot_s X^\bullet=\hom(\sk_s\Delta^\bullet,X^\bullet)\in\C$ where $\sk_s\Delta^\bullet\in\cS$ is the termwise $s$--skeleton of $\Delta^\bullet$.  Since $\Delta^\bullet$ is Reedy cofibrant and its skeletal inclusions are Reedy cofibrations, $\Tot X^\bullet$ is fibrant and $\{\Tot_s X^\bullet\} _{s\geq 0}$ is a tower of fibrations in $\C$ by 2.7.

For $M,Y\in\C$ and $n\geq 0$, let 
$$\pi_n(Y;M)~=~[M,Y]_n~=~[\Sigma^n M,Y]$$
denote the group or set of homotopy classes from $\Sigma^n M$ to $Y$ in the homotopy category $\Ho\C$.  Note that $\pi_n(Y;M)=\pi_n\map(\check{M},\bar{Y})$ where $\check{M}$ is a cofibrant replacement of $M$ and $\bar{Y}$ is a fibrant replacement of $Y$.

\subsection{The homotopy spectral sequence}
\label{sec:2.9}
As in \cite[pages 258 and 281]{BK}, the Tot tower $\{\Tot_s X^\bullet\} _{s\geq 0}$ now has a \emph{homotopy spectral sequence} 
$\{ E^{s,t}_r(X^\bullet;M)\}$ for $r\geq 1$ and $t\geq s\geq 0$, abutting to $\pi_{t-s}(\Tot X^\bullet;M)$ with differentials
$$d_r\co E^{s,t}_r(X^\bullet;M)\xrightarrow{\qquad}  E^{s+r,t+r-1}_r(X^\bullet;M)$$
and with natural isomorphisms
$$ E^{s,t}_1(X^\bullet;M)~\iso~\pi_{t-s}(\Fib_s X^\bullet;M)~\iso~N^s\pi_t(X^\bullet;M)$$
$$ E^{s,t}_2(X^\bullet;M)~\iso~\pi_{t-s}(\Fib_s X^\bullet;M)^{(1)}~\iso~\pi^s\pi_t(X^\bullet;M)$$
for $t\geq s\geq 0$ involving the fiber $\Fib_s X^\bullet$ of $\Tot_s X^\bullet\to\Tot_{s-1}X^\bullet$, the normalization 
$N^s(-)$, the couple derivation $(-)^{(1)}$, and the cosimplicial cohomotopy $\pi^s(-)$ (see \cite[2.2]{Bou1989} and \cite[page 284]{BK}).  This is equivalent to the ordinary homotopy spectral sequence of the cosimplicial space $\map(\check{M},X^\bullet)\in\cS_*$, and its basic properties follow immediately from earlier work.  We refer the reader to \cite[pages 261--264]{BK} and 
\cite[pages 63--67]{Bou1989}  for convergence results concerning the 
natural surjections
$\pi_i(\Tot X^\bullet;M)\to\lim_s Q_s\pi_i(\Tot X^\bullet;M)$ for $i\geq 0$ where $Q_s\pi_i(\Tot X^\bullet;M)$ denotes the image of $\pi_i(\Tot X^\bullet;M)\to\pi_i(\Tot_s X^\bullet;M)$ and concerning the natural inclusions
$E^{s,t}_{\infty+}(X^\bullet;M)\subset E^{s,t}_\infty(X^\bullet;M)$ where\break $E^{s,t}_{\infty+}(X^\bullet;M)$ denotes the kernel of  $Q_s\pi_{t-s}(\Tot X^\bullet;M)\to Q_{s-1}\pi_{t-s}(\Tot X^\bullet;M)$ and where \ 
$ E^{s,t}_{\infty}(X^\bullet;M) =  \bigcap_{r>s}  E^{s,t}_r(X^\bullet;M)$.  As in \cite{Bou1989}, the spectral sequence may be partially extended beyond the $t\geq s\geq 0$ sector, and there is an associated obstruction theory.  Finally, in preparation for our work on resolution model categories, we consider the following:

%old sec:2.10 omitted

\subsection{The external simplicial structure on $\cC$}
\label{sec:2.11}
For a category $\C$ with finite limits and colimiits, the category $\cC=s(\C^{\rm op})$ has an \emph{external simplicial structure} as in  Quillen \cite[II.1.7]{Qui} with a \emph{mapping space} $\map^c(X^\bullet,Y^\bullet)\in\Scal$, a \emph{cotensor} $\hom^c(K,X^\bullet)\in\cC$, and a \emph{tensor} $X^\bullet\tensor_c K\in\cC$ for $X^\bullet,Y^\bullet\in\cC$ and finite $K\in\Scal$. The latter are given  by 
$$\hom^c(K,X^\bullet)^n ~=~ \hom(K_n,X^n)$$ 
$$(X^\bullet\tensor_c K)^n~=~X^\bullet\tensor_{\Deltabf}(K\times\Delta^n)$$ 
for $n\geq 0$, using the coend over $\Deltabf$ and letting $\hom(S,X^n)$ and $X^n\tensor S$ respectively denote the product and coproduct of copies of $X^n$ indexed by a set $S$.  When $\C$ is a model category, the external simplicial structure on  $\cC$ will usually be incompatible with the Reedy model category structure.  However, it will satisfy the weakened version of {\bf SM7$^\prime$} obtained by replacing ``either $i$ or $j$ is trivial'' by ``$i$ is trivial'' (see \cite[page  372]{GJ}).  Moreover, as suggested by Meyer \cite[Theorem 2.4]{Mey}, we have the following:

\begin{lem}
\label{lem:2.12}
Suppose $\C$ is a bicomplete simplicial model category.  Then for $Y^\bullet\in\cC$ and $K\in\Scal$, there is a natural isomorphism
$$\Tot\hom^c(K,Y^\bullet)~\iso~\hom(K,\Tot Y^\bullet)\in\C.$$
\end{lem}

\begin{proof} 
It suffices adjointly to show, for $A\in\C$ and $K\in\Scal$, that there is a natural isomorphism 
$(A\tensor\Delta^\bullet)\tensor_c K \iso (A\tensor K)\tensor\Delta^\bullet\in \cC$.  This follows from the isomorphisms
$$(A\tensor\Delta^\bullet)\tensor_{\Deltabf} (K\times \Delta^n) ~\iso~ A\tensor(K\times\Delta^n)~\in ~\C$$
in codimensions $n\geq 0$, obtained by applying $A\tensor-$ to $\Delta^\bullet\tensor_{\Deltabf} (K\times\Delta^n)\iso K\times\Delta^n \in\Scal.$
\end{proof}

\subsection{The external homotopy relation}
\label{sec:2.13}
In a general simplicial category, two maps $f,g\co X\to Y$ are \emph{simplicially homotopic}  when $[f]=[g]$ in 
$\pi_0\map(X,Y)$.  In $cC$, to avoid ambiguity, we say that two maps $f,g\co X^\bullet\to Y^\bullet$ are \emph{externally homotopic} or \emph{cosimplicially homotopic} (written$f\overset{c}\sim g$) when $[f]=[g]$ in $\pi_0\map^c(X,Y)$.  For homomorphisms $\alpha,\beta\co A^\bullet\to B^\bullet$ of cosimplicial abelian groups, the relation $\alpha\overset{c}\sim\beta$ corresponds to the chain homotopy relation for $N\alpha,N\beta\co NA^\bullet\to NB^\bullet$ by Dold--Puppe \cite[Satz 3.31]{DP}, and hence $\alpha\overset{c}\sim\beta$ implies 
$\alpha_*=\beta_*\co \pi^sA^\bullet\to\pi^sB^\bullet$ for $s\geq 0$.  Likewise for homomorphisms $\alpha,\beta\co A^\bullet\to B^\bullet$ of cosimplicial groups (or pointed sets), the relation $\alpha\overset{c}\sim\beta$ implies 
$\alpha_*=\beta_*\co \pi^sA^\bullet\to\pi^sB^\bullet$ for $s=0,1$ (or $s=0$).

Over a bicomplete simplicial model category $\C$, we now have the following:

\begin{prop}
\label{prop:2.14}
If  $f,g\co X^\bullet\to Y^\bullet\in\cC$ \  are maps of Reedy fibrant objects with $f\overset{c} \sim g$, then \ 
$\Tot f, \Tot g\co \Tot X^\bullet\to\Tot Y^\bullet$ are simplicially homotopic.  Moreover, when $\C$ is pointed, \  
$f_*=g_*\co \pi_*(\Tot X^\bullet;M)\to\pi_*(\Tot Y^\bullet;M)$  \ and \ $f_*=g_*\co E^{s,t}_r(X^\bullet;M)\to E^{s,t}_r(Y^\bullet;M)$ for $M\in\C$, \  $t\geq s\geq 0$, and $2\leq r\leq \infty+$.
\end{prop}

\begin{proof}
$\Tot f$ and $\Tot g$ are simplicially homotopic since $\Tot$ preserves strict homotopies $X^\bullet\to\hom^c(\Delta^1,Y^\bullet)$  by Lemma \ref{lem:2.12}.  The proposition now follows by \ref{sec:2.13}.
\end{proof}

\section{Existence of resolution model categories}
\label{sec:3}

We now turn to the resolution model category structures of Dwyer--Kan--Stover \cite{DKS} on the category $\cC=s(\C^{op})$ of cosimplicial objects over a model category $\C$.  These have more weak equivalences than the Reedy structures and are much more flexible since they depend on a specified class of injective models in $\Ho C$.  Moreover, they are compatible with the external simplicial structure on $\cC$.  Our version of this theory is more general than the original one, and we have recast the proofs accordingly.  We must assume that our model category $\C$ is  \emph{left proper}, meaning that each pushout of a weak equivalence along a cofibration is a weak equivalence.  As explained in \cite[11.1]{Hir}, this condition holds for most familiar model categories including those whose objects are all cofibrant as assumed in \cite{DKS}.  For simplicity, we now also assume that $\C$ is pointed, and postpone the unpointed generalization until Section \ref{sec:12}.

\subsection{$\G$--injectives}
\label{sec:3.1}
Let $\C$ be a left proper pointed model category, and let $\G$ be a class of group objects in the homotopy category $\Ho\C$.  A map $i\co A\to B$ in $\Ho\C$ is called 
\emph{$\G$--monic} when $i^*\co [B,G]_n\to[A,G]_n$ is onto for each $G\in\G$ and $n\geq 0$, and an object $Y\in\Ho\C$ is called \emph{$\G$--injective} when $i^*\co [B,Y]_n\to[A,Y]_n$ is onto for each $\G$--monic map $i\co A\to B$ in $\Ho\C$ and $n\geq 0$.  For instance, the objects $\Omega^nG\in\Ho\C$ are $\G$--injective for $G\in\G$ and $n\geq 0$, and so are the retracts of their products.  The classes of $\G$--monic maps and of $\G$--injective objects in $\Ho\C$ clearly determine each other.  We say that $\Ho\C$ \emph{has enough $\G$--injectives} when each object of $\Ho\C$ is the source of a $\G$--monic map to a $\G$--injective target, and we then call $\G$ a \emph{class of injective models} in $\Ho\C$.  We always assume that a class of injective models consists of group objects in the homotopy category.  We say that an object of $\C$ is $\G$--\emph{injective} when it is $\G$--injective in $\Ho\C$, and say that a map in $\C$ is $\G$--\emph{monic} when it is $\G$--monic in $\Ho\C$.  In Lemma \ref{lem:3.7} below, we show that a fibrant object $F\in\C$ is $\G$--injective if and only if the fibration $F\to *$ has the right lifting property for the $\G$--monic cofibrations in $\C$.  Extending this condition, we say that a fibration in $\C$ is $\G$--\emph{injective} when it has the right lifting property for the $\G$--monic cofibrations in $\C$.  A more explicit characterization of $\G$--injective fibrations is given later in Lemma \ref{lem:3.10}.

\subsection{The $\G$--resolution model structure on $\cC$}
\label{sec:3.2}
Recall that a homomorphism in the category $sGrp$ of simplicial groups is a \emph{weak equivalence} or \emph{fibration} when its underlying map in $\Scal$ is one.  For a map $f\co X^\bullet\to Y^\bullet$ in $\cC$, we say:
\begin{enumerate}
\item[(i)]$f$ is a \emph{$\G$--equivalence} when $f^*\co [Y^\bullet,G]_n\to[X^\bullet,G]_n$ is a weak equivalence in $sGrp$ for each $G\in\G$ and $n\geq 0$;
\item[(ii)]$f$ is a \emph{$\G$--cofibration} when $f$ is a Reedy cofibration and $f^*\co [Y^\bullet,G]_n\to[X^\bullet,G]_n$ is a fibration in $sGrp$ for each $G\in\G$ and $n\geq 0$;
\item[(iii)]$f$ is a \emph{$\G$--fibration} when $f\co X^n\to Y^n\times_{M^nY^\bullet}M^nX^\bullet$ is a $\G$--injective fibration in $\C$ for $n\geq 0$.
\end{enumerate}
We let $cC^\G$ denote the category $\cC$  with weak equivalences defined as  $\G$--equivalences, with cofibrations defined as  $\G$--cofibrations, with fibrations defined as $\G$--fibrations, and with the external simplicial structure (\ref{sec:2.11}).

\begin{thm}[after Dwyer--Kan--Stover]
\label{thm:3.3}
If $\C$ is a left proper pointed model category with a class $\G$ of injective models in $\Ho\C$, then $\cC^\G$ is a left proper pointed simplicial model category.
\end{thm}

We call $\cC^\G$ the \emph{$\G$--resolution model category} and devote the rest of Section \ref{sec:3} to proving this theorem.   Since the proof is very long,  the reader might wish to proceed directly to Section \ref{sec:4} for a discussion of the result with some general examples.  We start by noting the following:

\begin{prop}
\label{prop:3.4}
The limit axiom {\bf MC1}, the weak equivalence axiom {\bf MC2}, and the retraction axiom {\bf MC3} hold in $cC^\G$.
\end{prop}

To go further, we must study $\G$--monic cofibrations and $\G$--injective fibrations in $\C$, and we start with a lemma due essentially to Dan Kan (see \cite[11.1.16]{Hir}).  It applies to a commutative diagram
$$
\begin{CD}
\Tilde{A}   @>u>>  A  @>>>  X\\
@VV\Tilde{i}V     @VViV     @VVpV\\
\Tilde{B}  @>v>>   B  @>>>  Y
\end{CD}
$$
in a left proper model category $\C$ such that $u$ and $v$ are weak equivalences, $\tilde{i}$ and $i$ are cofibrations, and $p$ is a fibration.

\begin{lem}
\label{lem:3.5}
If the combined square has a lifting $\tilde{B}\to X$, then the right square has a lifting $B\to X$.
\end{lem}
\begin{proof}
Using a lifting $\tilde{B}\to X$, we break the right square into 
$$
\begin{CD}
A @>>>  A\coprod_{\tilde{A}}\tilde{B}  @>>>  X\\
@VVV          @VVV                            @VVV\\
B @>>>  B                                 @>>>  Y
\end{CD}
$$
Since $\C$ is left proper, the maps  $\tilde{B} \to A\coprod_{\tilde{A}}\tilde{B} \to B$ are weak equivalences, and the second map factors into a trivial cofibration $A\coprod_{\tilde{A}}\tilde{B} \to E$ and trivial fibration $E\to B$.  Thus the original right square has a lifting $B\to E\to X$.
\end{proof}

Henceforth, we assume that $\C$ and $\G$ satisfy the hypotheses of Theorem \ref{thm:3.3}.  Since each cofibration $A\to B$ in $\C$ can be approximated by a cofibration $\check{A}\to\check{B}$ between cofibrant objects, Lemma \ref{lem:3.5} implies the following:

\begin{lem}
\label{lem:3.6}
A fibration in $\C$ is $\G$--injective if and only if it has the right lifting property for each $\G$--monic cofibration between cofibrant objects.
\end{lem}

This easily implies the following:

\begin{lem}
\label{lem:3.7}
A fibrant object $F\in\C$ is $\G$--injective in $\Ho\C$ if and only if the fibration $F\to *$ is $\G$--injective.
\end{lem}

The classes of $\G$--monic cofibrations and of $\G$--injective fibrant objects (or $\G$--injective fibrations) in $\C$ now determine each other by the following:

\begin{lem}
\label{lem:3.8}
A cofibration $i\co A\to B$ in $\C$ is $\G$--monic if and only if\break
$i^*\co \Hom_\C(B,F)\to \Hom_\C(A,F)$ is onto for each $\G$--injective fibrant object $F\in\C$.
\end{lem}
\begin{proof}
For the \emph{if} part, it suffices to show that $i^*\co [B,\Omega^nG]\to[A,\Omega^nG]$ is onto for each $G\in\G$ and $n\geq 0$.  Since $\C$ is left proper, each map $A\to\Omega^nG\in\Ho\C$ can be represented by a map $f\co A\to F\in\C$ for some $\G$--injective fibrant object $F$ with $F\homeq\Omega^nG$.  Since $f$ is in the image of $i^*\co \Hom_\C(B,F)\to\Hom_\C(A,F)$, the \emph{if} part follows easily, and the \emph{only if} part is trivial.
\end{proof}

\begin{lem}
\label{lem:3.9}
A map $f\co A\to B$ in $\C$ can be factored into a $\G$--monic cofibration $f^\prime\co A\to E$ and a $\G$--injective fibration 
$f^{\prime\prime}\co E\to B$.
\end{lem}
\begin{proof}
Since $\C$ is left proper and $\Ho\C$ has enough $\G$--injectives, we may choose a $\G$--monic cofibration 
$\gamma\co A\to F$ to a $\G$--injective fibrant object $ F\in\C$.  We factor $(f,\gamma)\co A\to B\times F$ as the composite of a cofibration $f^\prime\co A\to E$ and a trivial fibration $q\co E\to B\times F$.  This gives the desired factorization $f=f^{\prime\prime}f^\prime$ where $f^{\prime\prime}$ is the composite of $q$ with the projection $B\times F\to B$.
\end{proof}

As suggested by Paul Goerss, this leads to a fairly explicit characterization of the $\G$--injective fibrations in $\C$. A map $E\to Y$ in $\C$ is called \emph{$\G$--cofree} if it may be expressed as a composition of a trivial fibration $E\to Y\times F$ and a projection $Y\times F\to Y$ for some $\G$--injective fibrant object $F$.

\begin{lem}
\label{lem:3.10}
A map $X\to Y$ in $\C$ is a $\G$--injective fibration if and only if it is a retract of some $\G$--cofree map $E\to Y$.
\end{lem}

\begin{proof}
For the \emph{only if } part, we assume that $X\to Y$ is a $\G$--injective fibration, and we factor it as a composition of a $\G$--monic cofibration $X\to E$ and a $\G$--cofree map $E\to Y\times F\to Y$ as above.  Since $X\to Y$ has the right lifting property for the $\G$--monic cofibration $X\to E$, it must be a retract of the $\G$--cofree map $E\to Y$ as required.  This gives the \emph{only if} part, and the \emph{if} part is trivial.
\end{proof}

%old lemma 3.11 omitted

Finally, consider a push-out square in $\C$:  
$$
\begin{CD}
A    @>>>  C\\
@VViV     @VVjV\\
B    @>>>  D
\end{CD}
$$

\begin{lem}
\label{lem:3.12}
Suppose $i$ is a $\G$--monic cofibration in $\C$. Then so is $j$, and the functor $[-,G]_n$ carries the square to a pullback of groups for each $G\in\G$ and $n\geq 0$.
\end{lem}

\begin{proof}
The first conclusion follows by Lemma \ref{lem:3.8}, while the second follows homotopically since $\C$ is left proper and each $G\in\G$ is a group object in $\Ho\C$.
\end{proof}

Our next goal is to describe the $\G$--cofibrations of $\cC$ in terms of the $\G$--monic cofibrations of $\C$ using the following:

\subsection{Partial latching objects}
\label{sec:3.13}
For $X^\bullet\in\cC$ and a finite $K\in\Scal$, we obtain an object 
$X^\bullet\tensor_{\Deltabf}K =(X^\bullet\tensor_c K)^0\in\C$ as in \ref{sec:2.11}.  This gives the latching object
$L^nX^\bullet = X^\bullet\tensor_{\Deltabf}\partial\Delta^n$ as well as $X^n=X^\bullet\tensor_{\Deltabf}\Delta^n$ for $n\geq 0$.  We now let $L^n_kX^\bullet = X^\bullet\tensor_{\Deltabf}V^n_k$ for $n\geq k\geq 0$ where $V^n_k\subset\Delta^n$ is the $k$--horn spanned by $d_i\iota_n$ for all $i\ne k$. More generally, for a subset $\sigma\subset\{0,1,\cdots,n\}$, we let $L^n_\sigma X^\bullet=X^\bullet\tensor_{\Deltabf}F^n_\sigma$ where $F^n_{\sigma}\subset\Delta^n$ is spanned by $d_i\iota_n$ for all $i\in\sigma$.  Thus,  $L^n_k X^\bullet=L^n_\sigma X^\bullet$ for $\sigma = \{0,\dots,\hat{k},\dots,n\}$, although usually $L^n_k X^\bullet\neq L^n_{\{k\}} X^\bullet$. For a cofibration $J\to K$ of finite objects in $\Scal$ and a Reedy cofibration $X^\bullet\to Y^\bullet$ in $\cC$, we note that the map
$$(X^\bullet\tensor_{\Deltabf}K)\coprod_{X^\bullet\tensor_{\Deltabf}J}(Y^\bullet\tensor_{\Deltabf}J)\xrightarrow{\qquad} Y^\bullet\tensor_{\Deltabf}K$$
is a cofibration in $\C$ since
$$(X^\bullet\tensor_c K)\coprod_{X^\bullet\tensor_c J}(Y^\bullet\tensor_c J)\xrightarrow{\qquad} Y^\bullet\tensor_c K$$
is a Reedy cofibration in $\cC$.

\begin{prop}
\label{prop:3.14}
Let $f\co X^\bullet\to Y^\bullet$ be a Reedy cofibration in $\cC$.  Then:
\begin{enumerate}
\item[\rm(i)]$f$ is a $\G$--cofibration if and only if the cofibration $X^n\coprod_{L^n_k X^\bullet}L^n_k Y^\bullet\to Y^n$ is $\G$--monic whenever $n\geq k \geq 0$;
\item[\rm(ii)]$f$ is a $\G$--trivial cofibration if and only if the cofibration $X^n\coprod_{L^n X^\bullet}L^n Y^\bullet$\break $\to Y^n$ is $\G$--monic whenever $n\geq 0$.
\end{enumerate}
\end{prop}

\begin{proof}
For $G\in\G$, $\sigma\subset\{0,1,\cdots,n\}$, and $n\geq 0$, we obtain a square
$$
\begin{CD}
[Y^n,G]_*            @>{\rm Id}>>      [Y^n,G]_*\\
@VVV                             @VVV\\
[X^n\coprod_{L^n_\sigma X^\bullet}L^n_\sigma Y^\bullet,G]_*  @>>> 
             [X^n,G]_*\times_{M^\sigma_n [X^\bullet,G]_*}M^\sigma_n[Y^\bullet,G]_*
\end{CD}
$$
where $M^\sigma_n$ is the matching functor, dual to $L^n_\sigma$, for simplicial groups.  Each of the statements in ($i$) (resp.\ ($ii$)) asserts the surjectivity of a vertical arrow in this square for $\sigma$ of cardinality $|\sigma|=n$ (resp.\ 
$|\sigma|=n+1$).  The proposition now follows inductively using our next lemma.
\end{proof}

\begin{lem}
\label{lem:3.15}
Given $n\geq 1$, suppose that the cofibration $X^m\coprod_{L^m_\sigma X^\bullet} L^m_\sigma Y^\bullet\to Y^m$ is $\G$--monic for each $m<n$ and each $\sigma\subset\{0,1,\cdots,m\}$ with $|\sigma|=m$ (resp.\ $|\sigma|=m+1$).  Then the map
$$[X^m\textstyle\coprod_{L^m_\sigma X^\bullet}L^m_\sigma Y^\bullet,G]_*  \xrightarrow{\qquad}
             [X^m,G]_*\times_{M^\sigma_m [X^\bullet,G]_*}M^\sigma_m[Y^\bullet,G]_*$$
is an isomorphism for each $G\in\G$, each $m\leq n$, and each $\sigma\subset\{0,1,\cdots,m\}$ with $|\sigma|\leq m$
(resp.\ $|\sigma|\leq m+1$).
\end{lem}

\begin{proof}
We first claim that the cofibration $X^m\coprod_{L^m_\sigma X^\bullet} L^m_\sigma Y^\bullet\to Y^m$ is $\G$--monic for each $m<n$ and each $\sigma\subset\{0,1,\cdots,m\}$ with $|\sigma|\leq m$ (resp.\ $|\sigma|\leq m+1$).  This follows by inductively applying the first part of Lemma \ref{lem:3.12} to the pushout  squares
$$
\begin{CD}
X^{m-1}\coprod_{L^{m-1}_\sigma X^\bullet} L^{m-1}_\sigma Y^\bullet     @>>>
                            X^m\coprod_{L^m_\sigma X^\bullet} L^m_\sigma Y^\bullet\\
@VVV                                         @VVV\\
Y^{m-1}          @>>>         X^m\coprod_{L^m_\tau X^\bullet} L^m_\tau Y^\bullet
\end{CD}
$$
where $\sigma=\{i_1,\cdots,i_{k-1}\}$ and $\tau=\{i_1,\cdots,i_k\}$ for  $0\leq i_1<\cdots<i_k\leq m$ with $m<n$.  The lemma now follows  by inductively applying the pullback part of Lemma \ref{lem:3.12} to these squares with $m\leq n$.
\end{proof}

Proposition \ref{prop:3.14} combines with Lemma \ref{lem:3.8} to give the following:

\begin{cor}
\label{cor:3.16}
Let $f\co X^\bullet\to Y^\bullet$ be a Reedy cofibration in $\cC$.  Then $f$ is a $\G$--cofibration (resp.\ $\G$--trivial cofibration) if and only if  $f^*\co \Hom_{\C}(Y^\bullet,F)\to\Hom_{\C}(X^\bullet,F)$ is a fibration (resp.\ trivial fibration) in $\Scal$ for each $\G$--injective fibrant object $F\in\C$.                                                
\end{cor}

The $\G$--trivial cofibration condition on a map $X^\bullet\to Y^\bullet$ in $\cC$ now reduces to the $\G$--monic cofibration condition on each $X^n\coprod_{L^nX^\bullet}L^nY^\bullet\to Y^n$, just as  the $\G$--fibration condition reduces to the $\G$--injective fibration condition on each $X^n\to Y^n\times_{M^nY^\bullet}M^nX^\bullet$.  Hence the model category axioms                             pertaining to these conditions now follow easily.  

\begin{prop}
\label{prop:3.17}
The lifting and factorization axioms {\bf MC4}(ii) and {\bf MC5}(ii) (for fibrations and trivial cofibrations) hold in $\cC^{\G}$.
\end{prop}

\begin{proof}
This follows by Reedy's constructions \cite{Ree} since the $\G$--injective fibrations have the right lifting property for $\G$--monic cofibrations, and since the maps in $\C$ may be factored as in Lemma \ref{lem:3.9}.
\end{proof}

Using the external simplicial structure  (\ref{sec:2.11}) on $\cC$, we now also have the simplicial axiom {\bf SM7}$^\prime$  by the following:

\begin{prop}
\label{prop:3.18}
If $i\co A^\bullet\to B^\bullet\in \cC$ is a $\G$--cofibration and $j\co J\to K\in\Scal$ is a cofibration of finite objects, then the map
$$(A^\bullet\tensor_c K)\coprod_{A^\bullet\tensor_cJ}(B^\bullet\tensor_c J)\xrightarrow{\qquad} B^\bullet\tensor_cK$$
is a $\G$--cofibration in $\cC$ which is trivial if either $i$ or $j$ is trivial.
\end{prop}

\begin{proof}
Since this map is a Reedy cofibration by \ref{sec:2.11}, the result follows from Corollary \ref{cor:3.16} by an adjunction argument using the isomorphism $\Hom_{\C}(A^\bullet\tensor_c K,F) \iso \map(K, \Hom_{C}(A^\bullet,F))$ in $\Scal$ for $F\in\C$.
\end{proof}

To prove the factorizaton axiom {\bf MC5}(i) (for $\G$--cofibrations and $\G$--trivial fibrations), we need the following:

\begin{lem}
\label{lem:3.19}
The $\G$--cofibrations and $\G$--trivial cofibrations in $\cC$ are closed under pushouts.
\end{lem}

\begin{proof}This follows from Corollary \ref{cor:3.16}.
\end{proof}

Since the $\G$--cofibrant objects of $\cC$ are the same as the Reedy cofibrant ones, we may simply call them \emph{cofibrant}.

\begin{lem}
\label{lem:3.20}
A map $f\co X^\bullet\to Y^\bullet$ of cofibrant objects in $\cC$ can be factored into a $\G$--cofibration 
$i\co X^\bullet\to M_f$ and a $\G$--equivalence $q\co M_f\to Y^\bullet$.
\end{lem}

\begin{proof}
Let $M_f$ be the mapping cylinder
$$M_f\quad=\quad(X^\bullet\tensor_c\Delta^1){\textstyle\coprod_{X^\bullet}}Y^\bullet\quad =\quad (X^\bullet\tensor_c\Delta^1)\coprod_{X^\bullet\coprod X^\bullet}(Y^\bullet{\textstyle\coprod X^\bullet})$$
Then the natural map $i\co X^\bullet\to M_f$ is a $\G$--cofibration by Lemma \ref{lem:3.19} since $X^\bullet\coprod X^\bullet\to X^\bullet\tensor_c\Delta^1$ is a $\G$--cofibration by Proposition \ref{prop:3.18}.  Likewise, the natural map $j\co Y^\bullet\to M_f$ is a $\G$--trivial cofibration, and its natural left inverse $q\co M_f\to Y^\bullet$ is a $\G$--equivalence.  This gives the required factorization $f=qi$.
\end{proof}

We can now prove {\bf MC5}(i).

\begin{prop}
\label{prop:3.21}
A map $f\co X^\bullet\to Y^\bullet$ in $\cC$ can be factored into a $\G$--cofibration $i\co X^\bullet\to N_f$ and a $\G$--trivial fibration $p\co N_f\to Y^\bullet$.
\end{prop}

\begin{proof}
First take Reedy cofibrant replacements to give a map $\check{f}\co \check{X}^\bullet\to \check{Y}^\bullet$ and use Lemma \ref{lem:3.20} to factor $\check{f}$.  Then use a pushout of $\check{f}$ to factor $f$ into a $\G$--cofibration $j\co X^\bullet\to E^\bullet$ and a $\G$--equivalence $r\co E^\bullet\to Y^\bullet$.  Finally apply Proposition \ref{prop:3.17} to factor $r$ into a $\G$--trivial cofibration $s\co E^\bullet\to N_f$ and a $\G$--trivial fibration $p\co N_f\to Y^\bullet$, and let $i=sj$.  
\end{proof}

To prove the lifting axiom {\bf MC4}(i) (for $\G$--cofibrations and $\G$--trivial fibrations), we need several preliminary results.

\begin{lem}
\label{lem:3.22}
If a map $f$ in $\cC$ has the right lifting property for $\G$--cofibrations (resp.\ $\G$--trivial cofibrations), then $f$ is a $\G$--trivial fibration (resp.\ $\G$--fibration).
\end{lem}

\begin{proof}
This follows by first using Proposition \ref{prop:3.21} (resp.\ Proposition \ref{prop:3.17}) to factor $f$, and then using the given right lifting property to express $f$ as a retract of the appropriate factor.
\end{proof}

\begin{lem}
\label{lem:3.23}
For a $\G$--fibrant object $F^\bullet\in\cC$ and a cofibration (resp.\ trivial cofibration)  $L\to K$ of finite objects in $\Scal$, the induced map $\hom^c(K,F^\bullet)\to\hom^c(L,F^\bullet)\in cC$ has the right lifting property for $\G$--trivial cofibrations (resp.\ $\G$--cofibrations).
\end{lem}

\begin{proof}
This follows by Propositions \ref{prop:3.17} and  \ref{prop:3.18}.  
\end{proof}

We now let $PF^\bullet\in\cC$ be the standard \emph{path object} given by
$$PF^\bullet\quad=\quad\hom^c(\Delta^1,F^\bullet)~\times_{F^\bullet}~*\quad=\quad\hom^c(\Delta^1,F^\bullet)~
\times_{F^\bullet\times F^\bullet}~F^\bullet.$$

\begin{lem}
\label{lem:3.24}
For a $\G$--fibrant object $F^\bullet\in\cC$, the natural map $PF^\bullet\to F^\bullet$ (resp.\ $PF^\bullet\to *$) has the right lifting property for $\G$--trivial cofibrations (resp.\ $\G$--cofibrations) in $\cC$.
\end{lem}

\begin{proof}
This follows from Lemma \ref{lem:3.23} since right lifting properties are preserved by pullbacks.
\end{proof}

\begin{lem}
\label{lem:3.25}
If $F^\bullet\to *$ is a $\G$--trivial fibration with $F^\bullet$ cofibrant, then  $F^\bullet\to *$ has the right lifting property for $\G$--cofibrations.
\end{lem}

\begin{proof}
The $\G$--fibration $PF^\bullet\to F^\bullet$ has a cross-section by Proposition \ref{prop:3.17}, and $F^\bullet\to *$ has the right lifting property for $\G$--cofibrations since $PF^\bullet\to *$ does by Lemma \ref{lem:3.24}.  
\end{proof}

\begin{lem}
\label{lem:3.26}
If $f$ and $g$ are maps in $\cC$ such that $gf$ is a $\G$--cofibration and $f$ is a Reedy cofibration, then $f$ is a $\G$--cofibration.
\end{lem}

\begin{proof}
This follows since a simplicial group homomorphism $G\to H$ is a fibration if and only if it induces surjections of Moore normalizations $N_qG\to N_qH$ for $q>0$ (see \cite[II$\S$3]{Qui}).
\end{proof}

We can now prove {\bf MC4}(i).

\begin{prop}
\label{prop:3.27}
A $\G$--trivial fibration $f\co X^\bullet\to Y^\bullet$ in $\cC$ has the right lifting property for $\G$--cofibrations.
\end{prop}

\begin{proof}
First suppose that $X^\bullet$ is cofibrant.  By Proposition \ref{prop:3.21}, the map $X^\bullet\to*$ factors into a $\G$--cofibration $\phi\co X^\bullet\to F^\bullet$ and a $\G$--trivial fibration $F^\bullet\to*$, and the map $(f,\phi)\co X^\bullet\to Y^\bullet\times F^\bullet$ factors into a Reedy cofibration $X^\bullet\to E^\bullet$ and a Reedy trivial fibration $E^\bullet\to Y^\bullet\times F^\bullet$.  Then the map $E^\bullet\to Y^\bullet$ is a $\G$--trivial fibration with the right lifting property for $\G$--cofibrations by Lemmas \ref{lem:3.22} and \ref{lem:3.25}.  Hence, $X^\bullet\to E^\bullet$ is a $\G$--equivalence and a $\G$--cofibration by Lemma \ref{lem:3.26}.  Thus $X^\bullet\to Y^\bullet$ is a retract of $E^\bullet\to Y^\bullet$ by Proposition \ref{prop:3.17}, and $X^\bullet\to Y^\bullet$ inherits the right lifting property for $\G$--cofibrations.  In general, by Lemma \ref{lem:3.5} (applied in Reedy's $\cC$), it suffices to show that $X^\bullet\to Y^\bullet$ has the right lifting property for each $\G$--cofibration of cofibrant objects $C^\bullet\to D^\bullet$.  This follows since a map $C^\bullet\to X^\bullet$ factors into a Reedy cofibration $C^\bullet\to\check{X}^\bullet$ and Reedy trivial fibration $\check{X}^\bullet\to X^\bullet$, where the composed map $\check{X}^\bullet\to X^\bullet\to Y^\bullet$ must have the right lifting property for $\G$--cofibrations since it is a $\G$--trivial fibration with $\check{X}^\bullet$ cofibrant.
\end{proof}

This completes the proof that $\cC^{\G}$ is a simplicial model category, and Theorem \ref{thm:3.3} will follow from the following:

\begin{prop}
\label{prop:3.28}
The $\G$--resolution model category $\cC^{\G}$ is left proper.
\end{prop}

\begin{proof}
By \cite[Lemma 9.4]{Bou2000}, it suffices to show that a pushout of a $\G$--equivalence $f\co A^\bullet\to Y^\bullet$ along a $\G$--cofibration $A^\bullet\to B^\bullet$ of cofibrant objects is a $\G$--equivalence.  We may factor $f$ into a $\G$--equivalence $\phi\co A^\bullet\to\check{Y}^\bullet$ with $\check{Y}^\bullet$ cofibrant and a Reedy weak equivalence $q\co \check{Y}^\bullet\to Y^\bullet$.  The proposition now follows since the pushout of $\phi$ is a $\G$--equivalence by \cite[Theorem B]{Ree}, and the pushout of $q$ is a Reedy weak equivalence.
\end{proof}

\section{Examples of resolution model categories}
\label{sec:4}

If $\C$ is a left proper pointed model category with a class $\G$ of injective models in $\Ho\C$, then Theorem \ref{thm:3.3} gives the $\G$--resolution model category $\cC^{\G}$.  In this section, we discuss some general examples of these model categories.

\subsection{Dependence of $\cC^{\G}$ on $\G$}
\label{sec:4.1}
As initially defined, the $\G$--resolution model structure on $\cC$ seems to depend strongly on $\G$.  However, by Proposition \ref{prop:3.14}, the $\G$--cofibrations and $\G$--trivial cofibrations in $\cC$ are actually determined by the $\G$--monic maps in $\Ho\C$.  Hence, \emph{the $\G$--resolution model structure on $\cC$ is determined by the class of $\G$--monic maps, or equivalently by the class of $\G$--injective objects in $\Ho\C$.}  

\subsection{A refinement of Theorem \ref{thm:3.3}}
\label{sec:4.2}
Adding to the hypotheses of Theorem \ref{thm:3.3}, we suppose that the model category $\C$ is factored (\ref{sec:2.1}) and that the class $\G$ of injective models is \emph{functorial}, meaning that there exists a functor $\Gamma\co \C\to\C$ and a transformation $\gamma\co 1_{\C}\to\Gamma(X)$ such that $\gamma\co X\to\Gamma(X)$ is a $\G$--monic map to a $\G$--injective object $\Gamma(X)$ for each $X\in\C$.  Then the model category $\cC^{\G}$ is also factored by the constructions in our proof of Theorem \ref{thm:3.3}.  Of course, if $\C$ is bicomplete, then $\cC^{\G}$ is also bicomplete.

\subsection{Constructing $\cC^{\G}$ for discrete $\C$}
\label{sec:4.3}
Let $\C$ be a pointed category with finite limits and colimits, and give $\C$ the \emph{discrete} model category structure in which the weak equivalences are the isomorphisms, and the cofibrations and fibrations are arbitrary maps.  Then $\Ho\C=\C$ with $[X,Y]_0=\Hom_{\C}(X,Y)$ and with $[X,Y]_n=*$ for $X,Y\in\C$ and $n>0$.  Now let $\G$ be a class of group objects in $\C$.  If $\C$ has enough $\G$--injectives, then we have  a simplicial model category $\cC^{\G}$ by Theorem \ref{thm:3.3}.  This provides a dualized variant of Quillen's Theorem 4 in \cite[II$\S$4]{Qui}, allowing many possible choices of ``relative injectives'' in addition to Quillen's canonical choice.  For instance, we consider the following:

\subsection{Abelian examples}
\label{sec:4.4}
Let $\C$ be an abelian category, viewed as a discrete model category, and let $\G$ be a class of objects in $\C$ such that $\C$ has enough $\G$--injectives. Recall that $\cC$ is equivalent to the category $Ch^+\C$ of nonnegatively graded cochain complexes over $\C$ by the Dold--Kan correspondence (see eg \cite{DP} or \cite{GJ}).  Thus the $\G$--resolution model category $\cC^{\G}$ corresponds to a model category $Ch^+\C^{\G}$. For a cochain map $f\co X\to Y$ in $Ch^+\C^{\G}$, a careful analysis shows that:
\begin{enumerate}
\item[(i)]$f$ is a \emph{$\G$--equivalence} when $f^*\co H_n\Hom(Y,G)\iso H_n\Hom(X,G)$ for each $G\in\G$ and $n\geq 0$;
\item[(ii)]$f$ is a \emph{$\G$--cofibration} when $f\co X^n\to Y^n$ is $\G$--monic for $n\geq 1$;
\item[(iii)]$f$ is a \emph{$\G$--fibration} when $f\co X^n\to Y^n$ is splittably epic with a $\G$--injective kernel for $n\geq 0$.
\end{enumerate}
For example, when $\C$ has enough injectives and $\G$ consists of them all, we recover Quillen's model category $Ch^+\C^{\G}$ \cite[II$\S$4]{Qui} where: (i) the $\G$--equivalences are the cohomology equivalences; (ii) the $\G$--cofibrations are the maps monic in positive degrees; and (iii) the $\G$--fibrations are the epic maps with injective kernels in all degrees.  For another example, when $\G$ consists of all objects in $\C$, we obtain a model category $Ch^+\C^{\G}$ where: (i) the $\G$--equivalences are the chain homotopy equivalences; (ii) the $\G$--cofibrations are the maps splittably monic in positive degrees; and (iii) the $\G$--fibrations are the maps splittably epic in all degrees.  In this example, all cochain complexes are $\G$--fibrant and $\G$--cofibrant.

\subsection{Constructing $\cC^{\G}$ for small $\G$}
\label{sec:4.5}
Let $\C$ be a left proper pointed model category with arbitrary products, and let $\G$ be a (small) set of group objects in $\Ho\C$.  Then $\Ho\C$ has enough $\G$--injectives, since for each $X\in \Ho\C$, there is a natural $\G$--monic map
$$X\xrightarrow{\qquad} \prod_{G\in\G}\quad\prod_{n\geq 0}~\prod_{f\co X\to\Omega^n G}\Omega^n G$$ 
to a $\G$--injective target, where $f$ ranges over all maps $X\to\Omega^nG$ in $\Ho\C$.  Thus we have a simplicial model category $\cC^{\G}$ by Theorem \ref{thm:3.3}.  Note that an object $X\in\Ho\C$ is $\G$--injective if and only if $X$ is a retract of a product of terms $\Omega^nG$ for various $G\in\G$ and $n\geq 0$.  Also note that if $\C$ is factored, then the class $\G$ is functorial by a refinement of the above construction, and hence the model category $\cC^{\G}$ is factored by \ref{sec:4.2}.

\subsection{A homotopical example}
\label{sec:4.6}
Let $\Ho_*=\Ho\Scal_*$ be the pointed homotopy category of spaces, and recall that a cohomology theory $E^*$ is representable by spaces $\underline{E}_n\in\Ho_*$ with $\tilde{E}^nX\iso[X,\underline{E}_n]$ for $X\in\Ho_*$ and $n\in\Z$.  
For $\G=\{\underline{E}_n\}_{n\in\Z}$, we obtain a $\G$--resolution model category $\cS_*^{\G}$ by \ref{sec:4.5}.  Note that the $\G$--equivalences in $\cS_*$ are the maps inducing $\pi_s\tilde{E}^*$--isomorphisms for $s\geq 0$.  Also note that $\cS_*^{\G}$ is factored by \ref{sec:4.5}.  Our next example will involve the following:

\subsection{Quillen adjoints}
\label{sec:4.7}
Let $\C$ and $\D$ be left proper pointed model categories, and let $S\co \C\leftrightarrows\D\thinspace\colon T$ be \emph{Quillen adjoint} functors, meaning that $S$ is left adjoint to $T$ and the following  equivalent conditions are satisfied:
(i) $S$ preserves cofibrations and $T$ preserves fibrations;
(ii) $S$ preserves cofibrations and trivial cofibrations; and
(iii) $T$ preserves fibrations and trivial fibrations.
Then by \cite{Qui} or \cite[Theorem 9.7]{DS}, $S$ has a \emph{total left derived functor} $\Lcal S\co \Ho\C\to\Ho\D$, and $T$ has a \emph{total right derived functor} $\R T\co \Ho\D\to\Ho\C$, where $\Lcal S$ is left adjoint to $\R T$.  Moreover, 
$\Lcal S$ preserves homotopy cofiber sequences and suspensions, while $\R T$ preserves homotopy fiber sequences and loopings.  

\subsection{Construction $\cC^{\G}$ from Quillen adjoints}
\label{sec:4.8}
Let $S\co\C\leftrightarrows\D\,\colon T$ be Quillen adjoints as in \ref{sec:4.7}, and let $\Hcal$ be a class of injective models in $\Ho\D$.  Then we obtain a class $\G=\{(\R T)H~|~H\in\Hcal\}$ of injective models in $\Ho C$, and obtain Quillen adjoints $S\co\cC^{\G}\leftrightarrows\cD^{\Hcal}\,\colon T$.  We note that if $\C$ and $\D$ are factored and $\Hcal$ is functorial, then $\G$ is also functorial and hence $\cC^{\G}$ and $\cD^{\Hcal}$ are factored.

\subsection{Another homotopical example}
\label{sec:4.9}
Let $\Scal p$ be the model category of spectra in the sense of \cite{BF} (see also \cite{HSS}), and let $\Ho^s=\Ho(\Scal p)$ be the stable homotopy category.  The infinite suspension and 0-space functors $\Scal_*\leftrightarrows\Scal p$ are Quillen adjoints, and their total derived functors are the usual infinite suspension and infinite loop functors
$\Sigma^{\infty}\co \Ho_*\leftrightarrows\Ho^s\,\colon\Omega^{\infty}$.  Let $S\in\Ho^s$ be the sphere spectrum, and suppose that $E\in\Ho^s$ is a \emph{ring spectrum}, meaning that it is equipped with a multiplication $E\wedge E\to E\in\Ho^s$ and unit $S\to E\in\Ho^s$ satisfying the identity and associativity properties in $\Ho^s$.  Let $\Hcal$ be the class of $E$--module spectra in $\Ho^s$ and note that $\Ho^s$ has enough $\Hcal$--injectives since the unit maps $X\to E\wedge X$ are $\Hcal$--monic with $\Hcal$--injective targets.  Thus by \ref{sec:4.8}, we obtain a class $\G=\{\Omega^{\infty}N~|~N\in\Hcal\}$ of injective models in $\Ho_*$, and we have resolution model categories $\cS p^{\Hcal}$ and $\cS_*^{\G}$  by Theorem \ref{thm:3.3}.  Various alternative choices of $\G$ will lead to the same $\G$--injectives in $\Ho_*$ and hence to the same resolution model category $\cS_*^{\G}$.  For instance, we could equivalently let $\G$ be $\{\Omega^{\infty}(E\wedge\Sigma^{\infty}X)~|~X\in\Ho_*\}$ or $\{\Omega^{\infty}(E\wedge Y)~|~Y\in\Ho^s\}$.  These resolution model categories are factored.

\section{Derived functors, completions and homotopy spectral sequences}
\label{sec:5}

Let $\C$ be a left proper pointed model category with a class $\G$ of injective models in $\Ho\C$.  We now introduce $\G$--resolutions of objects in $\C$ and use them to construct right derived functors, completions, and the associated homotopy spectral sequences.  In Section \ref{sec:6}, we shall see that a weaker sort of $\G$--resolution will suffice for these purposes.

\subsection{$\G$--resolutions in $\C$}
\label{sec:5.1}
A \emph{$\G$--resolution} (= \emph{cosimplicial $\G$--injective resolution}) of an object $A\in\C$ consists of a $\G$--trivial cofibration $\alpha\co A\to\bar{A}^\bullet$ to a $\G$--fibrant object $\bar{A}^\bullet$ in $\cC$, where $A$ is considered constant in $\cC$.  This exists for each $A\in\C$ by {\bf MC5} in $\cC^{\G}$, and exists functorially when $\cC^{\G}$ is factored.  In general, $\G$--resolutions are natural up to external homotopy (\ref{sec:2.13}) by the following:

\begin{lem}
\label{lem:5.2}
If $\alpha\co A\to I^\bullet$ is a $\G$--trivial cofibration  in $\cC^{\G}$, and if $f\co A\to J^\bullet$ is a map to a $\G$--fibrant object $J^\bullet\in\cC^{\G}$, then there exists a map $\phi\co I^\bullet\to J^\bullet$ with $\phi\alpha=f$ and $\phi$ is unique up to external homotopy.
\end{lem}

\begin{proof}
This follows since $\alpha^*\co \map^c(I^\bullet,J^\bullet)\to\map^c(A,J^{\bullet})$ is a trivial fibration in $\Scal$ by {\bf SM7} in $\cC^{\G}$.
\end{proof}

The terms of a $\G$--resolution are $\G$--injective by the following:

\begin{lem}
\label{lem:5.3}
If an object $I^\bullet\in\cC$ is $\G$--fibrant, then $I^n$ is $\G$--injective and fibrant in $\C$ for $n\geq 0$.
\end{lem}

\begin{proof}
More generally, if $f\co X^\bullet\to Y^\bullet$ is a $\G$--fibration in $\cC$, then  
$f\co X^n\to Y^n\times_{M^nY^\bullet}M^nX^\bullet$  is a $\G$--injective fibration for $n\geq 0$ by definition, and hence each $f\co X^n\to Y^n$ is a $\G$--injective fibration in $\C$ by Corollary 2.6 of \cite[page 366]{GJ}.
\end{proof}

Conseqently, the terms $I^n$ are H-spaces in $\Ho\C$ by the following: 

\begin{lem}
\label{lem:5.4}
If $J$ is a $\G$--injective object in $\Ho\C$, then $J$ admits a multiplication with unit.
\end{lem}

\begin{proof}
The coproduct-to-product map $J\vee J\to J\times J$ is $\G$--monic since $\Omega^nG$ is a group object of $\Ho\C$ for each $G\in\G$ and $n\geq 0$.  Hence, the folding map $J\vee J\to J$ extends to a map $J\times J\to J$ giving the desired multiplication.
\end{proof}

\subsection{Right derived functors}
\label{sec:5.5}
Let $T\co \C\to\M$ be a functor to an abelian category $\M$.  We define the  \emph{right derived functor} 
$\R^s_{\G}T\co \C\to\M$ for $s\geq 0$, with a natural transformation $\epsilon\co T\to\R^0_{\G}T$, by setting
$\R^s_{\G}T(A)~=~\pi^s T\bar{A}^\bullet~=~H^s(NT\bar{A}^\bullet)$
for $A\in\C$, where $A\to\bar{A}^\bullet\in\cC$ is a $\G$--resolution of $A$ and $NT\bar{A}^\bullet$ is the normalized cochain complex of $T\bar{A}^\bullet\in c\M$.  This is well-defined up to natural equivalence by \ref{sec:2.13} and \ref{lem:5.2}.  Similarly, let $U\co \C\to Grp$ and $V\co \C\to Set_*$ be functors to the categories of groups and pointed sets.  We define the \emph{right derived functors} $\R^0_{\G}U\co \C\to Grp$ and $\R^1_{\G}U,\R^0_{\G}V\co \C\to Set_*$ by setting $\R^s_{\G}U(A)=\pi^sU\bar{A}^\bullet$ and $\R^s_{\G}V(A)=\pi^s V\bar{A}^\bullet$ as above.  Since the  $\G$--fibrant objects in $\cC$ are termwise $\G$--injective by Lemma \ref{lem:5.3}, these derived functors 
depend only on the restrictions of $T$, $U$, $V$ to the full subcategory of $\G$--injective objects in $\C$.  Thus they may be defined for such restricted functors.  

\subsection{Abelian examples}
\label{sec:5.6}
Building on \ref{sec:4.4}, suppose $\C$ is an abelian category with a class $\G$ of injective models, and suppose $T\co \C\to\M$ is a functor to an abelian category $\M$.  Then a $\G$--resolution of $A\in\C$ corresponds to an augmented cochain complex $A\to\tilde{A}^\bullet\in Ch^+\C$ where $\tilde{A}^n$ is $\G$--injective for $n\geq 0$ and where the augmented chain complex $\Hom(\tilde{A}^\bullet,G)\to\Hom(A,G)$ is acyclic for each $G\in\G$.  When $T$ is additive, we have $\R^s_{\G}T(A)=H^sT\tilde{A}^\bullet$ for $s\geq 0$, and we recover the usual right derived  functors 
$\R^s_{\G}T\co \C\to\M$ of relative homological algebra \cite{EM}. In general, we obtain relative versions of the Dold--Puppe \cite{DP} derived functors.

Now suppose that the model category $\C$ is simplicial and bicomplete.

\subsection{$\G$--completions}
\label{sec:5.7}
For an object $A\in\C$, we define the \emph{$\G$--completion} 
$\alpha:A\to\hat{L}_{\G}A\in\Ho\C$ by setting $\hat{L}_{\G}A=\Tot\bar{A}^\bullet$ where $A\to\bar{A}^\bullet\in\cC$ is a $\G$--resolution of $A$.  This determines a functor $\hat{L}_{\G}\co \C\to\Ho\C$ which is well-defined up to natural equivalence by \ref{lem:5.2} and \ref{prop:2.14}.  In fact, by Corollary \ref{cor:8.2} below, the $\G$--completion will give a functor $\hat{L}_{\G}\co \Ho\C\to\Ho\C$ and a natural transformation $\alpha\co {\rm Id}\to\hat{L}_{\G}$ belonging to a triple on $\Ho\C$.  When $\C$ is factored and $\G$ is functorial (\ref{sec:4.2}), the $\G$--completion is canonically represented by a functor $\hat{L}_{\G}\co \C\to\C$ with a natural transformation $\alpha\co {\rm Id}\to\hat{L}_{\G}$.

\subsection{$\G$--homotopy spectral sequences} 
\label{sec:5.8}
For objects $A,M\in\C$, we define the \emph{$\G$--homotopy spectral sequence} 
$$\{E^{s,t}_r(A;M)_{\G}\}_{r\geq 2}$$ of $A$ with coefficients $M$ by setting $E^{s,t}_r(A;M)_{\G}=E^{s,t}_r(\bar{A}^\bullet;M)$ for $0\leq s \leq t$ and $2\leq r\leq\infty+$ using the homotopy spectral sequence (\ref{sec:2.9}) of $\bar{A}^\bullet$ for a $\G$--resolution $A\to\bar{A}^\bullet$.  Since this is the homotopy spectral sequence of a pointed cosimplicial space \ $\map(\check{M},\bar{A}^\bullet)$, composed of H-spaces by \ref{lem:5.4}, we see that 
$E^{s,t}_r(A;M)_{\G}$ is a pointed set for $0\leq s=t\leq r-2$ and is otherwise an abelian group by \cite[Section 2.5]{Bou1989}.  The spectral sequence is fringed on the line $t=s$ as in \cite{BK}, and the differentials
$$d_r\co E^{s,t}_r(A;M)_{\G}\xrightarrow{\qquad} E^{s+r,t+r-1}_r(A;M)_{\G}$$
are homomorphisms for $t>s$.  It has 
$$E^{s,t}_2(A;M)_{\G}~=~\pi^s\pi_t(\bar{A}^\bullet;M)~=~\R^s_{\G}\pi_t(A;M)$$ 
for $0\leq s\leq t$ by \ref{sec:2.9} and \ref{sec:5.5}, and it abuts to $\pi_{t-s}(\hat{L}_{\G}A;M)$ with the usual convergence properties which may be expressed using the natural surjections
$\pi_i(\hat{L}_{\G}A;M)\to\lim_s Q_s\pi_i(\hat{L}_{\G}A;M)$ for $i\geq 0$ and the natural inclusions
$$E^{s,t}_{\infty+}(A;M)_{\G}\subset E^{s,t}_\infty(A;M)_{\G}$$
as in 2.9.  The spectral sequence is well-defined up to natural equivalence and depends functorially on $A,M\in\C$ by  \ref{lem:5.2} and  \ref{prop:2.14}.

\subsection{Immediate generalizations}
\label{sec:5.9}
The above notions extend to an arbitrary object $A^\bullet\in\cC$ in place of $A\in\C$.  A \emph{$\G$--resolution} of $A^\bullet$ still consists of a $\G$--trivial cofibration $\alpha\co A\to\bar{A}^\bullet$ to a $\G$--fibrant object $\bar{A}^\bullet\in\cC$.  A functor $T\co \C\to\M$ to an abelian category $\M$ still has \emph{right derived functors} $\R^s_{\G}T\co \cC\to\M$ with $\R^s_{\G}T(A^\bullet)=\pi^sT\bar{A}^\bullet\in\M$ for $s\geq 0$.  Moreover, $A^\bullet$ still has a \emph{$\G$--homotopy spectral sequence} $\{E^{s,t}_r(A^\bullet;M)_{\G}\}_{r\geq 2}$ with coefficients $M\in\C$, where 
$E^{s,t}_r(A^\bullet;M)_{\G}=E^{s,t}_r(\bar{A}^\bullet;M)$ for $0\leq s\leq t$ and $2\leq r\leq\infty+$.  This has 
$$E^{s,t}_2(A^\bullet;M)_{\G}~=~\pi^s\pi_t(\bar{A}^\bullet;M)~=~\R^s_{\G}\pi_t(A^\bullet;M)$$
for $t\geq  s\geq 0$ and abuts to $\pi_{t-s}\Tot_{\G}A^\bullet$ where $\Tot_{\G}A^\bullet=\Tot\bar{A}^\bullet\in\Ho\C$ (see \ref{prop:8.1}).  It retains the properties described above in \ref{sec:5.8}.

\section{Weak resolutions are sufficient}
\label{sec:6}

Let $\C$ be a left proper pointed model category with a class $\G$ of injective models in $\Ho\C$.  We now introduce the weak $\G$--resolutions of objects in $\C$ and show that they may be used in place of actual $\G$--resolutions to construct right derived functors, $\G$--completions, and $\G$--homotopy spectral sequences.  This is convenient since weak $\G$--resolutions arise naturally from triples on $\C$ (see Section \ref{sec:7}) and are generally easy to recognize.

\begin{defn}
\label{defn:6.1}
A \emph{weak $\G$--resolution} of an object $A\in\C$ consists of a $\G$--equivalence $A\to Y^\bullet$ in $\cC$ such that $Y^n$ is $\G$--injective for $n\geq 0$.  Such a $Y^\bullet$ is called \emph{termwise $\G$--injective}.
\end{defn}

Any $\G$--fibrant object of $\cC$ is termwise $\G$--injective by Lemma \ref{lem:5.3}, and hence any $\G$--resolution is a weak $\G$--resolution.  As our first application, we consider the right derived functors of a functor $T\co \C\to\N$ where $\N$ is an abelian category or $\N=Grp$ or $\N=Set_*$.  We suppose that $T$ carries weak equivalences in $\C$ to isomorphisms in $\N$.

\begin{thm}
\label{thm:6.2}
If $A\to Y^\bullet\in\cC$ is a weak $\G$--resolution of an object $A\in\C$, then there is a natural isomorphism $\R^s_{\G}T(A)\iso\pi^sTY^\bullet$ for $s\geq 0$.
\end{thm}

It is understood that  $s=0,1$ when $\N=Grp$ and that $s=0$ when $\N=Set_*$. This theorem will be proved in \ref{sec:6.14}, and we cite two elementary consequences.

\begin{cor}
\label{cor:6.3}
If $A\in\C$ is $\G$--injective, then $\epsilon\co T(A)\iso\R^0_{\G}T(A)$ and $\R^s_{\G}T(A)$\break $=0$ for $s>0$.
\end{cor}

\begin{proof}
This follows using the weak $\G$--resolution ${\rm Id}\co A\to A$.
\end{proof}

A map $f\co A\to B$ in $\C$ is called a \emph{$\G$--equivalence} when $f^*\co [B,G]_n\iso[A,G]_n$ for $G\in\G$ and $n\geq 0$, or equivalently when $f$ is a $\G$--equivalence of constant objects in $\cC$.

\begin{cor}
\label{cor:6.4}
If $f\co A\to B$ is a $\G$--equivalence in $\C$, then $f_*\co \R^s_{\G}T(A)\iso\R^s_{\G}T(B)$ for $s\geq 0$.
\end{cor}

\begin{proof}
This follows since $f$ composes with a weak $\G$--resolution of $B$ to give a weak $\G$--resolution of $A$.
\end{proof}

To give similar results for $\G$--completions and $\G$--homotopy spectral sequences, we suppose that $\C$ is simplicial and bicomplete.

\begin{thm}
\label{thm:6.5}
Suppose $A\to Y^\bullet$ is a weak $\G$--resolution of an object $A\in\C$. Then there is a natural equivalence 
$\hat{L}_{\G}A\homeq\Tot\underline{Y}^\bullet\in\Ho\C$ for a Reedy fibrant replacement $\underline{Y}^\bullet$ of $Y^\bullet$, and  there are natural isomorphisms  $E^{s,t}_r(A;M)_{\G}\iso E^{s,t}_r(\underline{Y}^\bullet;M)$ and $Q_s\pi_i(\hat{L}_{\G}A;M)\iso Q_s\pi_i(\Tot\underline{Y}^\bullet;M)$ for $M\in C$,  $0\leq s\leq t$,\ $2\leq r\leq\infty+$, and $i\geq 0$.
\end{thm}

This will be proved later in \ref{sec:6.19} and partially generalized in \ref{thm:9.5}.  It has the following elementary consequences.

\begin{cor}
\label{cor:6.6}
Suppose $A\in\C$ is $\G$--injective. Then $\hat{L}_{\G}A\homeq A\in\Ho\C$ and 
$$
E^{s,t}_r(A;M)_{\G}~\iso~
\begin{cases}
\pi_t(A;M)        &\text{when $s=0$}\\
0                 &\text{when $0<s\leq t$}
\end{cases}
$$
for $M\in\C$  and  $2\leq r\leq\infty+$.
\end{cor}

\begin{cor}
\label{cor:6.7}
If $f\co A\to B$ is a $\G$--equivalence in $\C$, then $f$ induces $\hat{L}_{\G}A\homeq \hat{L}_{\G}B$ and 
$E^{s,t}_r(A;M)_{\G}\iso E^{s,t}_r(B;M)_{\G}$ for $M\in\C$, $0\leq s\leq t$, and $2\leq r\leq\infty+$.
\end{cor}

In particular, the $\G$--completion $\hat{L}_{\G}\co \C\to\Ho\C$ carries weak equivalences to equivalences and induces a functor $\hat{L}_{\G}\co \Ho\C\to\Ho\C$.  To prepare for the proofs of Theorems \ref{thm:6.2} and \ref{thm:6.5}. we need the following:

\subsection{The model category $c(\cC^{\G})$}
\label{sec:6.8}
Let $c(\cC^{\G})$ be the Reedy model category of cosimplicial objects $X^{\bullet\bullet}=\{X^{n\bullet}\}_{n\geq 0}$ over the $\G$--resolution model category $\cC^{\G}$.  Its structural maps are called \emph{Reedy $\G$--equivalences}, \emph{Reedy $\G$--cofibrations}, and \emph{Reedy $\G$--fibrations}.  Thus a map $f\co X^{\bullet\bullet}\to Y^{\bullet\bullet}$ is a Reedy $\G$--equivalence if and only if $f\co X^{n\bullet}\to Y^{n\bullet}$ is a $\G$--equivalence for each $n\geq 0$.  Moreover, if $f\co X^{\bullet\bullet}\to Y^{\bullet\bullet}$ is a Reedy $\G$--cofibration (resp.\ Reedy $\G$--fibration), then $f\co X^{n\bullet}\to Y^{n\bullet}$ is a $\G$--cofibration (resp.\ $\G$--fibration) for each $n\geq 0$ by \cite[Corollary VII.2.6]{GJ}.  Let $\diag\co c(\cC^{\G})\to\cC^{\G}$ be the functor with $\diag Y^{\bullet\bullet}=\{Y^{nn}\}_{n\geq 0}$.

\begin{lem}
\label{lem:6.9}
If $f\co X^{\bullet\bullet}\to Y^{\bullet\bullet}$ is a Reedy $\G$--equivalence, then $\diag f\co \diag X^{\bullet\bullet}$\break $\to\diag Y^{\bullet\bullet}$ is a $\G$--equivalence.
\end{lem}

\begin{proof}
For each $G\in\G$, the bisimplicial group hommorphism $f^*\co [Y^{\bullet\bullet},G]_*\to[X^{\bullet\bullet},G]_*$ restricts to a weak equivalence $[Y^{n\bullet},G]_*\to[X^{n\bullet},G]_*$ for $n\geq 0$, and thus restricts to a weak equivalence 
$[\diag Y^{\bullet\bullet},G]_*\to[\diag X^{\bullet\bullet},G]_*$ by \cite[Theorem B.2]{BF}.
\end{proof}

\begin{lem}
\label{lem:6.10}
If $f\co X^{\bullet\bullet}\to Y^{\bullet\bullet}$ is a Reedy ${\G}$--fibration, then $\diag f\co \diag X^{\bullet\bullet}$\break $\to\diag Y^{\bullet\bullet}$ is a $\G$--fibration.
\end{lem}

\begin{proof}
For $X^{\bullet\bullet}\in c(\cC^{\G})$, we may express $\diag X^{\bullet\bullet}$ as an end
$$
\diag X^{\bullet\bullet}~\iso~\int_{[n]\in\Deltabf}\hom^c(\Delta^n,Y^{n\bullet}),
$$
and hence interpret $\diag X^{\bullet\bullet}$ as the total object (\ref{sec:2.8}) of the cosimplicial object $X^{\bullet\bullet}$ over $\cC^{\G}$.  The lemma now follows  by \ref{sec:2.7}.
\end{proof}

\subsection{Special $\G$--fibrant replacements}
\label{sec:6.11}
For an object $Y^\bullet\in\cC^{\G}$, we let $\con Y^\bullet\in c(\cC^{\G})$ be the vertically constant object with 
$(\con Y^\bullet)^{n,i}=Y^n$ for $n,i\geq 0$.  We choose a Reedy $\G$--trivial cofibration 
$\alpha\co \con Y^\bullet\to\vec{Y}^{\bullet\bullet}$ to a Reedy $\G$--fibrant target $\vec{Y}^{\bullet\bullet}$, and we let 
$\Bar{\Bar{Y}}^\bullet=\diag\vec{Y}^{\bullet\bullet}$.  This induces a $\G$--equivalence $\alpha\co Y^\bullet\to\Bar{\Bar{Y}}^\bullet$ with $\Bar{\Bar{Y}}^\bullet$ $\G$--fibrant by Lemmas \ref{lem:6.9} and \ref{lem:6.10}.  With some work, we can show that this special $\G$--fibrant replacement $\alpha\co Y^\bullet\to\Bar{\Bar{Y}}^\bullet$ is actually a $\G$--resolution, but that will not be needed.

Let $T\co \C\to\M$ be a functor to an abelian category $\M$ such that $T$ carries weak equivalences to isomorphisms.

\begin{lem}
\label{lem:6.12}
If $Y^\bullet\in\cC^{\G}$ is Reedy fibrant and termwise $\G$--injective, then the above map $\alpha\co Y^\bullet\to\Bar{\Bar{Y}}^\bullet$ induces an isomorphism $\alpha_*\co \pi^sTY^\bullet\iso\pi^sT\Bar{\Bar{Y}}^\bullet$ for $s\geq 0$.
\end{lem}

\begin{proof}
Since $\alpha\co Y^n\to\vec{Y}^{n\bullet}$ is a $\G$--resolution of the $\G$--injective fibrant object $Y^n$, we have $\pi^sT\vec{Y}^{n\bullet}=0$ for $s>0$ and $\pi^0T\vec{Y}^{n\bullet}\iso TY^n$.  Hence, 
$T\alpha\co T(\con Y^\bullet)\to T\vec{Y}^{\bullet\bullet}$ restricts to $\pi^*$--equivalences of all vertical complexes, and must therefore restrict to a $\pi^*$--equivalence of the diagonal complexes by the Eilenberg--Zilber--Cartier theorem of Dold--Puppe \cite{DP}.  Hence, $\alpha\co T(Y^\bullet)\to T(\Bar{\Bar{Y}}^\bullet)$ is a $\pi^*$--equivalence.
\end{proof}

\begin{lem}
\label{lem:6.13}
If $Y^\bullet,Z^\bullet\in\cC^{\G}$ are termwise $\G$--injective and $f\co Y^\bullet\to Z^\bullet$ is a $\G$--equivalence, then $f_*\co \pi^sTY^\bullet\iso\pi^sTZ^\bullet$ for $s\geq 0$.
\end{lem}

\begin{proof}
After replacements, we may assume that $Y^\bullet$ and $Z^\bullet$ are Reedy fibrant.  Let $\alpha\co Y^\bullet\to\Bar{\Bar{Y}}^\bullet$ and $\beta\co Z^\bullet\to\Bar{\Bar{Z}}^\bullet$ be special $\G$--fibrant replacements as in \ref{sec:6.11} with an induced map $\Bar{\Bar{f}}\co \Bar{\Bar{Y}}^\bullet\to\Bar{\Bar{Z}}^\bullet$ such that $\Bar{\Bar{f}}\alpha=\beta f$.  Then $\alpha$ and $\beta$ are $\pi^*T$--equivalences by Lemma \ref{lem:6.12}.  After Reedy cofibrant replacements, $\Bar{\Bar{f}}$ becomes a $\G$--equivalence of $\G$--fibrant cofibrant objects and hence a cosimplicial homotopy equivalence.  Thus $\Bar{\Bar{f}}$ is also a $\pi^*T$--equivalence, and hence so is $f$.
\end{proof}

\subsection{Proof of Theorem \ref{thm:6.2}}
\label{sec:6.14}
Consider the case of \ $T\co \C\to\M$ as above.  Given a weak $\G$--resolution $\alpha\co A\to Y^\bullet$, we choose $\G$--resolutions $u\co A\to\bar{A}^\bullet$ and $v\co Y^\bullet\to\bar{Y}^\bullet$, and choose $\bar{\alpha}\co \bar{A}^\bullet\to\bar{Y}^\bullet$ with $\bar{\alpha}u=v\alpha$.  Then
$$
\R^s_{\G}TA~\iso~\pi^sT\bar{A}^\bullet~\iso~\pi^sT\bar{Y}^\bullet~\iso~\pi^sTY^\bullet
$$
for $s\geq 0$ by Lemma \ref{lem:6.13} as required.  The remaining cases of $T\co \C\to Grp$ and $T\co \C\to Set_*$ are similarly proved.\endproof

To prepare for the proof of Theorem \ref{thm:6.5}, we let $\M$ be a bicomplete simplicial model category.  For an object $M^{\bullet\bullet}\in c(\cM)$, we define $\Tot^vM^{\bullet\bullet}\in\cM$ by $(\Tot^vM^{\bullet\bullet})^n=\Tot(M^{n\bullet})$ for $n\geq 0$.

\begin{lem}
\label{lem:6.15}
For $M^{\bullet\bullet}\in c(\cM)$, there is a natural isomorphism
$\Tot\Tot^vM^{\bullet\bullet}$\break $\iso\Tot\diag M^{\bullet\bullet}$.
\end{lem}

\begin{proof}
The functor $\Tot\co \cM\to\M$ preserves inverse limits and gives\break $\Tot\hom^c(K,N^\bullet)\iso\hom(K,\Tot N^\bullet)$ for $N^\bullet\in\cM$ and $K\in\Scal$ by Lemma \ref{lem:2.12}.  Hence, the induced functor $\Tot^v=c(\Tot)\co c(\cM)\to\cM$ respects total objects (\ref{sec:2.8}), and we have
$$
\Tot\Tot^vM^{\bullet\bullet}~\iso~\Tot\Tot M^{\bullet\bullet}~\iso~\Tot\diag M^{\bullet\bullet}
$$
with $\Tot M^{\bullet\bullet}\iso\diag M^{\bullet\bullet}$ by the proof  of Lemma \ref{lem:6.10}.
\end{proof}

Using the Reedy and Reedy--Reedy model category structures (\ref{thm:2.3}) on $\cM$ and $c(\cM)$, we have the following:

\begin{lem}
\label{lem:6.16}
If $M^{\bullet\bullet}\to N^{\bullet\bullet}$ is a Reedy--Reedy fibration in $c(\cM)$, then 
$\Tot^v M^{\bullet\bullet}\to\Tot^vN^{\bullet\bullet}$ is a Reedy fibration in $\cM$.
\end{lem}

\begin{proof}
This follows since $\Tot\co \cM\to\M$ preserves fibrations and inverse limits.
\end{proof}

For Theorem \ref{thm:6.5}, we also need the following comparison lemma of \cite[6.3 and 14.4]{Bou1989} whose hypotheses are expressed using notation from $\S\S2,14$ of that paper.

\begin{lem}
\label{lem:6.17}
Let $f\co V^\bullet\to W^\bullet$ be a map of pointed fibrant cosimplicial spaces such that:
\begin{enumerate}
\item[\rm(i)]$f_*\co \pi^0\pi_0V^\bullet\iso\pi^0\pi_0W^\bullet$;
\item[\rm(ii)]$f$ induces an equivalence $\Tot\pi^{gd}_1V^\bullet\iso\Tot\pi^{gd}_1W^\bullet$ of groupoids;
\item[\rm(iii)]$f_*\co \pi^*\pi_t(V^\bullet,b)\iso\pi^*\pi_t(W^\bullet,fb)$ for each vertex $b\in\Tot_2V^\bullet$ and $t\geq 2$.
\end{enumerate}
Then $f$ induces an equivalence $\Tot V^\bullet\iso\Tot W^\bullet$ and isomorphisms 
$Q_s\pi_i\Tot V^\bullet$\break $\iso Q_s\pi_i\Tot W^\bullet$ and $E^{s,t}_rV^\bullet\iso E^{s,t}_rW^\bullet$ for $0\leq s\leq t$, $2\leq r\leq\infty+$, and $i\geq 0$.
\end{lem}

This leads to our final preparatory lemma.

\begin{lem}
\label{lem:6.18}
If $Y^\bullet\in\cC^{\G}$ is Reedy fibrant and termwise $\G$--injective,  then $\alpha\co Y^\bullet\to\Bar{\Bar{Y}}^\bullet$ (as in \ref{sec:6.11}) induces an equivalence $\Tot Y^\bullet\homeq \Tot\Bar{\Bar{Y}}^\bullet$ and isomorphisms $Q_s\pi_i(\Tot Y^\bullet;M)\iso Q_s\pi_i(\Tot\Bar{\Bar{Y}}^\bullet;M)$ and
$E^{s,t}_r(Y^\bullet;M)\iso E^{s,t}_r(\Bar{\Bar{Y}}^\bullet;M)$ for a cofibrant $M\in\C$,\ $0\leq s\leq t$,\  $2\leq r\leq\infty+$, and $i\geq 0$.
\end{lem}

\begin{proof}
Since $Y^n$ is $\G$--fibrant, the $\G$--resolution $\alpha\co Y^n\to\vec{Y}^{n\bullet}$ is a cosimplicial homotopy equivalence such that $Y^n$ is a strong deformation retract of $\vec{Y}^{n\bullet}$ for $n\geq 0$.  Thus $\alpha\co Y^\bullet\to\Tot^v\vec{Y}^{\bullet\bullet}$ is a Reedy weak equivalence of Reedy fibrant objects by Proposition \ref{prop:2.14} and Lemma \ref{lem:6.16}, and 
$$
\Tot Y^\bullet\longrightarrow\Tot\Tot^v\vec{Y}^{\bullet\bullet}~\iso~\Tot\Bar{\Bar{Y}}^\bullet
$$
is an equivalence by Lemma \ref{lem:6.15} as desired.  For the remaining conclusions, it suffices to show that $\map(M,Y^\bullet)\to\map(M,\Bar{\Bar{Y}}^\bullet)$ satisfies the hypotheses (i)--(iii) of Lemma \ref{lem:6.17}.  This follows by double complex arguments since $\map(M,Y^n)\to\map(M,\vec{Y}^{n\bullet})$ is a cosimplicial homotopy equivalence such that $\map(M,Y^n)$ is a strong deformation retract of $\map(M,\vec{Y}^{n\bullet})$ for $n\geq 0$, and hence this homotopy equivalence induces: (i) a $\pi^0\pi_0$--isomorphism; (ii) a $\Tot\pi_1^{gd}$--equivalence; and (iii) a 
$\pi^*\pi_t(-,b)$--isomorphism for each vertex $b\in\Tot_2\map(M,Y^\bullet)$ and $t\geq 2$.  In (iii) we note that the vertex $b$ determines a map $b\co \sk_2\Delta^n\to\map(M,Y^n)$ which provides a sufficiently well defined basepoint for $\map(M,Y^n)$ since the space $\sk_2\Delta^n$ is simply connected.
\end{proof}

\subsection{Proof of Theorem \ref{thm:6.5}}
\label{sec:6.19}
The proof of Theorem \ref{thm:6.2}  is easily adapted to give Theorem \ref{thm:6.5} using Lemma \ref{lem:6.18} in place of Lemma \ref{lem:6.12}. \endproof

\subsection{Immediate generalizations}
\label{sec:6.20}
In \ref{sec:5.9}, we explained how the notions of \emph{$\G$--resolution}, \emph{right derived functor}, and \emph{$\G$--homotopy spectral sequence} apply not merely to objects $A\in\C$ but also to objects $A^\bullet\in\cC$.  Similarly, we may now define a \emph{weak $\G$--resolution} of an object $A^\bullet\in\cC$ to be a \emph{$\G$--equivalence} $A^\bullet\to Y^\bullet$ such that $Y^\bullet$ is termwise $\G$--injective.  Then the results \ref{thm:6.2}--\ref{cor:6.7} have immediate generalizations where: $A,B\in\cC$ are replaced by $A^\bullet,B^\bullet\in\cC$; \ \emph{$\G$--injective} is replaced by \emph{termwise $\G$--injective}; and $\hat{L}_{\G}A$ is replaced by $\Tot_{\G}A^\bullet$.

\section{Triples give weak resolutions}
\label{sec:7}

We now explain how weak $\G$--resolutions may be constructed from suitable triples, and give some examples.  We can often show that our weak $\G$--resolutions are actual $\G$--resolutions, but that seems quite unnecessary.

\subsection{Triples and triple resolutions}
\label{sec:7.1}
Recall that a \emph{triple} or \emph{monad} $\langle\Gamma,\eta,\phi\rangle$ on a category $\M$ consists of a  functor $\Gamma\co \M\to\M$ with transformations $\eta\co 1_{\M}\to\Gamma$ and $\mu\co \Gamma\,\Gamma\to\Gamma$ satisfying the identity and associativity conditions.  For an object $M\in\M$, the \emph{triple resolution} $\alpha\co M\to\Gamma^\bullet M\in\cM$ is the augmented cosimplicial object with $(\Gamma^\bullet M)^n=\Gamma^{n+1}M$ and
\begin{align*}
d^i~&=~\Gamma^i\eta\Gamma^{n-i+1}\co (\Gamma^\bullet M)^n\to(\Gamma^\bullet M)^{n+1}\\
s^i~&=~\Gamma^i\mu\Gamma^{n-i}\co (\Gamma^\bullet M)^{n+1}\to(\Gamma^\bullet M)^n
\end{align*}
for $n\geq -1$.  The \emph{augmentation map} $\alpha\co M\to\Gamma^\bullet M\in\cM$ is given by $d^0\co M\to(\Gamma^\bullet M)^0$.  An object $I\in\M$ is called \emph{$\Gamma$--injective} if $\eta\co I\to\Gamma I$ has a left inverse.

\begin{lem}
\label{lem:7.2}
For a triple $\langle\Gamma,\eta,\mu\rangle$ on $\M$ and object $M\in\M$, the triple resolution $\alpha\co M\to\Gamma^\bullet M$ induces a weak equivalence $\alpha^*\co \Hom(\Gamma^\bullet M,I)\to\Hom(M,I)$ in $\Scal$ for each $\Gamma$--injective $I\in\M$.
\end{lem}

\begin{proof}
Since $I$ is a retract of $\Gamma I$, it suffices to show that $\alpha^*\co \Hom(\Gamma^\bullet M,\Gamma I)\to\Hom(M,\Gamma I)$ is a weak equivalence.  This follows by Lemma \ref{lem:7.3} below since the augmented simplicial set \ $\Hom(\Gamma^\bullet M,\Gamma I)$ admits a left contraction $s_{-1}$ with $s_{-1}f=\mu(\Gamma f)$ for each simplex $f$.
\end{proof}

For an augmented simplicial set $K$ with augmentation operator $d_0\co K_0\to K_{-1}$, a \emph{left contraction} consists of functions $s_{-1}\co K_n\to K_{n+1}$ for $n\geq -1$ such that, in all degrees, there are identities $d_0s_{-1}=1$, $d_{i+1}s_{-1}=s_{-1}d_i$ for $i\geq 0$, and $s_{j+1}s_{-1}=s_{-1}s_j$ for $j\geq-1$.  As shown in  \cite[page 190]{GJ}, we have the following:

\begin{lem}
\label{lem:7.3}
If $K$ admits a left contraction, then the augmentation map $K\to K_{-1}$ is a weak equivalence in $\Scal$.
\end{lem}

Now suppose that $\C$ is a left proper pointed model category with a given class $\G$ of injective models in $\Ho\C$.

\begin{thm}
\label{thm:7.4}
Let $\langle\Gamma,\eta,\mu\rangle$ be a triple on $\C$ such that $\eta\co A\to\Gamma A$ is $\G$--monic with $\Gamma A$ \  $\G$--injective for each $A\in\C$.  If $\Gamma\co \C\to\C$ preserves weak equivalences, then the triple resolution $\co A\to\Gamma^\bullet A$ is a weak $\G$--resolution for each $A\in\C$.
\end{thm}

\begin{proof}
Since $(\Gamma^\bullet A)^n=\Gamma^{n+1}A$ is  $\G$--injective for $n\geq 0$, it suffices to show that $\alpha\co A\to\Gamma^\bullet A$ induces a weak equivalence $\alpha^*\co [\Gamma^\bullet A,\Omega^tG]\to[A,\Omega^tG]$ in $\Scal$ for each $G\in\G$ and $t\geq 0$.  This follows by Lemma \ref{lem:7.2} since $\langle\Gamma,\eta,\mu\rangle$ gives a triple on $\Ho\C$ such that each $\Omega^t G$ is $\Gamma$--injective.
\end{proof}

Various authors including Barr--Beck \cite{BB}, 
Bousfield--Kan \cite{BK}, and Bendersky--Thompson \cite{BT} have used triple resolutions to define right derived functors, completions, or homotopy spectral sequences, and we can now fit these constructions into our framework.  Starting with a triple, we shall find a compatible class of injective models giving the following:

\subsection{An interpretation of triple resolutions}
\label{sec:7.5}
Let $\M$ be a left proper pointed model category, and let $\langle\Gamma,\eta,\mu\rangle$ be a triple on $\M$ such that $\Gamma$ preserves weak equivalences.  Then there is an induced triple on $\Ho\M$ which is also denoted by 
$\langle\Gamma,\eta,\mu\rangle$.  For each $X\in\Ho\M$, we suppose:
\begin{enumerate}
\item[(i)]$\Gamma X$ is a group object in $\Ho\M$; 
\item[(ii)]$\Omega\Gamma X$ is $\Gamma$--injective in $\Ho\M$.
\end{enumerate}
Now $\G=\{\Gamma X~|~X\in\Ho\M\}$ is a class of injective models in $\Ho\M$, and we can interpret the triple resolution $\alpha\co A\to\Gamma^\bullet A$ of $A\in\M$ as a weak $\G$--resolution by Theorem \ref{thm:7.4}.

\subsection{The discrete case}
\label{sec:7.6}
Suppose $\M$ is a pointed category with finite limits and colimits, and suppose $\langle\Gamma,\eta,\mu\rangle$ is a triple on $\M$ such that $\Gamma X$ is a group object in $\M$ for each $X\in\M$.  The above discussion now applies to the discrete model category $\M$  and allows us to interpret the triple resolution $\alpha\co A\to\Gamma^\bullet A$ of $A\in\M$ as a weak $\G$--resolution where $\G=\{\Gamma X~|~X\in\M\}$.  Thus if $T\co \C\to\N$ is a functor to an abelian category $\N$ or to $\N=Grp$ or to $\N=Set_*$, then we obtain $\R^s_{\G}T(A)=\pi^sT(\Gamma^\bullet A)$ thereby recovering the right derived functors of Barr-Beck \cite{BB} and others.

\subsection{The Bousfield--Kan resolutions}
\label{sec:7.7} For a ring $R$, there is a triple $\langle R,\eta,\mu\rangle$ on the model category $\Scal_*$ of pointed spaces where $(RX)_n$ is the free $R$--module on $X_n$ modulo the relation $[*]=0$.  This satisfies the conditions of \ref{sec:7.5}, so that we may interpret the Bousfield-Kan resolution $\alpha\co A\to R^\bullet A\in\cS_*$ as a weak $\G$--resolution of $A\in\Scal_*$ where $\G=\{RX~|~X\in\Ho_*\}$ or equivalently $\G=\{\Omega^\infty N~|~N\text{ is an $HR$--module spectrum}\}$ as in \ref{sec:4.9}.  Thus we recover the Bousfield-Kan $R$--completion $R_\infty X\homeq\hat{L}_{\G}X$ and the accompanying homotopy spectral sequence.  More generally, we consider

\subsection{The Bendersky--Thompson resolutions}  
\label{sec:7.8}
For a ring spectrum $E$, there is an obvious triple on $\Ho_*$ carrying a space $X$ to $\Omega^\infty(E\wedge\Sigma^\infty X)$.  In \cite[Proposition 2.4]{BT}, Bendersky and Thompson suppose that $E$ is represented by an $S$--algebra \cite{EKMM}, and they deduce that the above homotopical triple is represented by a topological triple, and hence by a triple $\langle E,\eta,\mu\rangle$ on $\Scal_*$.  This triple satisfies the conditions of \ref{sec:7.5}, so that we may interpret the Bendersky--Thompson resolution $A\to E^\bullet A\in\cS_*$ as a weak $\G$--resolution of $A\in\Scal_*$, where 
$\G=\{EX~|~X\in\Ho_*\}$ or equivalently (see \ref{sec:4.9}) where $\G$ is the class 
$\{\Omega^\infty N~|~N \text{ is an $E$--module spectrum}\}$ or the class $\{\Omega^\infty(E\wedge Y)~|~Y\in\Ho^s\}$.  Thus we recover the Bendersky--Thompson $E$--completion $\hat{X}_E\homeq\hat{L}_{\G}X$ and the accompanying homotopy spectral sequence $\{E^{s,t}_r(A;M)_E\}=\{E^{s,t}_r(A;M)_{\G}\}$ over an arbitrary ring spectrum $E$ which need not be an $S$--algebra. As pointed out by Dror Farjoun \cite[page 36]{Dro}, Libman \cite{Lib}, and Bendersky--Hunton \cite{BH}, this generality can also be achieved by using restricted cosimplicial $E$--resolutions without codegeneracies.  However, we believe that codegeneracies remain valuable; for instance, they are essential for our constructions of pairings and products in these spectral sequences \cite{Bou2002}.  We remark that these various alternative constructions of homotopy spectral sequences over a ring spectrum all produce equivalent $E_2$--terms  and almost surely produce equivalent spectral sequences from that level onward.  Finally we consider the following:

\subsection{The loop-suspension resolutions}
\label{sec:7.9}
For a fixed integer $n\geq 1$, we let $\langle\Gamma, \eta, \mu\rangle$ be a triple on $\Scal_*$ representing the $n$-th loop-suspension triple $\Omega^n\Sigma^n$ on $\Ho_*$.  This satisfies the conditions of \ref{sec:7.5}, so that we may interpret the $n$-th loop-suspension resolution $A\to\Gamma^\bullet A\in \cS_*$ as a weak $\G$--resolution of $A\in\Scal_*$ where $\G=\{\Omega^n\Sigma^nX~|~X\in\Ho_*\}$ or equivalently where $\G = \{\Omega^nY~|~Y\in\Ho_*\}$.  The \emph{n-{\rm th} loop-suspension completion} of A is now given by $\Tot\underline{\Gamma^\bullet A} \homeq \hat{L}_{\G}A$, and will be identified in \ref{sec:9.9}.

\section{Triple structures of completions}
\label{sec:8}

Let $\C$ be a left proper, bicomplete, pointed simplicial model category with a class $\G$ of injective models in $\Ho\C$.  We now show that the $\G$--completion functor $\hat{L}_{\G}\co \Ho\C\to\Ho\C$ and transformation $\alpha\co 1\to\hat{L}_{\G}$ belong to a triple on $\Ho\C$, and we introduce notions of $\G$--completeness, $\G$--goodness, and $\G$--badness in $\Ho\C$.  This generalizes work of Bousfield--Kan \cite{BK} on the $R$--completion functor $R_\infty\co \Ho_*\to\Ho_*$ where $R$ is a ring.

By \ref{sec:2.7} and \ref{sec:2.8}, the functor $\Delta^\bullet\tensor-\co \C\to\cC$ is left adjoint to $\Tot\co \cC\to\C$, and these functors become Quillen adjoint (\ref{sec:4.7}) when $\cC$ is given the Reedy model category structure.  This remains true when $\cC$ is given the $\G$--resolution model category structure by the following:

\begin{prop}
\label{prop:8.1}
The functors $\Delta^\bullet\tensor-\co \C\to\cC^{\G}$ and $\Tot\co \cC^{\G}\to\C$ are Quillen adjoint.
\end{prop}

\begin{proof}
For a cofibration (resp.\ trivial cofibration) $A\to B$ in $\C$, it suffices by Corollary \ref{cor:3.16} to show that the Reedy cofibration $A\tensor\Delta^\bullet\to B\tensor\Delta^\bullet$ induces a fibration (resp.\ trivial fibration) 
$\Hom_{\C}(B\tensor\Delta^\bullet,F)\to\Hom_{\C}(A\tensor\Delta^\bullet,F)$ in $\Scal$ for each $\G$--injective fibrant object $F\in\C$.  This follows from the axiom {\bf SM7} on $\C$, since this fibration is just \ $\map(B,F)\to\map(A,F)$.      
\end{proof}

The resulting adjoint functors
$$
\Lcal(-\tensor\Delta^\bullet)\co\Ho\C~\leftrightarrows~\Ho(\cC^{\G})\,\colon\R\Tot
$$
will be denoted by
$$
\con\co\Ho\C~\leftrightarrows~\Ho(\cC^{\G})\,\colon\Tot_{\G}.
$$
Thus, for $A\in\Ho\C$ and $X^\bullet\in\Ho(\cC^{\G})$, we have $\con(A)\homeq A\in\Ho(\cC^{\G})$ and $\Tot_{\G}X^\bullet\homeq\Tot\bar{X}^\bullet$ where $X^\bullet\to\bar{X}^\bullet$ is a $\G$--fibrant approximation to $X^\bullet$.

\begin{cor}
\label{cor:8.2}
The $\G$--completion functor $\hat{L}_{\G}\co \Ho\C\to\Ho\C$ and transformation $\alpha\co 1\to\hat{L}_{\G}$ belong to  a triple $\langle\hat{L}_{\G},\alpha,\mu\rangle$ on $\Ho\C$.
\end{cor}

\begin{proof}
We easily check that $\hat{L}_{\G}$ and $\alpha$ belong to the adjunction triple of the above functors \ $\con$ 
and \ $\Tot_{\G}$.
\end{proof}

\begin{defn}
\label{defn:8.3}
An object $A\in\Ho\C$ is called \emph{$\G$--complete} if $\alpha\co A\homeq\hat{L}_{\G}A$; $A$ is called \emph{$\G$--good} if $\hat{L}_{\G}A$ is $\G$--complete; and $A$ is called \emph{$\G$--bad} if $\hat{L}_{\G}A$ is not $\G$--complete.
\end{defn}

A $\G$--injective object of $\Ho\C$ is $\G$--complete by Corollary \ref{cor:6.6}, and a $\G$--complete object is clearly $\G$--good.  To study these properties, we need the following: 

\begin{lem}
\label{lem:8.4}
For a map $f\co A\to B$ in $\Ho\C$, the following are equivalent:
\begin{enumerate}
\item[\rm(i)]$f\co A\to B$ is a $\G$--equivalence (see \ref{cor:6.4});
\item[\rm(ii)]$\hat{L}_{\G}f\co \hat{L}_{\G}A\homeq\hat{L}_{\G}B$;
\item[\rm(iii)]$f^*\co [B,I]\iso[A,I]$ for each $\G$--complete object $I\in\Ho\C$.
\end{enumerate}
\end{lem}

\begin{proof}
We have $(i)\Rightarrow(ii)$ by Corollary \ref{cor:6.7}.  To show $(ii)\Rightarrow(iii)$, note that a map $u\co A\to I$ extends to a map 
$$
\alpha^{-1}_I(\hat{L}u)(\hat{L}f)^{-1}\alpha_B\co B\xrightarrow{\qquad}I
$$
so $f^*$ is onto; and note that if $u\co A\to I$ extends to a map $r\co B\to I$, then 
$$
r~=~\alpha^{-1}_I(\hat{L}r)\alpha_B~=~\alpha^{-1}_I(\hat{L}r)(\hat{L}f)(\hat{L}f)^{-1}\alpha_B~=~\alpha^{-1}_I(\hat{L}u)(\hat{L}f)^{-1}\alpha_B
$$
so $f^*$ is monic.  To show $(iii)\Rightarrow(i)$, note that $\Omega^nG$ is $\G$--complete for each $G\in\G$ and $n\geq 0$, since it is $\G$--injective.
\end{proof}

\begin{prop}
\label{prop:8.5}
An object $A\in\Ho\C$ is $\G$--good if and only if $\alpha\co A\to\hat{L}_{\G}A$ is a $\G$--equivalence.
\end{prop}

\begin{proof}
If either of the maps $\alpha,\hat{L}\alpha\co \hat{L}A\to\hat{L}\hat{L}A$ is an equivalence in $\Ho\C$, then so is the other since they have the same left inverse $\mu\co \hat{L}\hat{L}A\to\hat{L}A$.  The result now follows from Lemma \ref{lem:8.4}.
\end{proof}

Thus the $\G$--completion $\alpha\co A\to\hat{L}_{\G}A$ of a $\G$--good object $A\in\Ho\C$ may be interpreted as the localization of $A$ with respect to the $\G$--equivalences (see \cite[2.1]{Bou1975}), and the $\G$--completion functor is a reflector from the category of $\G$--good objects to that of $\G$--complete objects in $\Ho\C$.  In contrast, for $\G$--bad objects, we have the following:

\begin{prop}
\label{prop:8.6}
If an object $A\in\Ho\C$ is $\G$--bad, then so is $\hat{L}_{\G}A$.
\end{prop}

\begin{proof}
Using the triple structure $\langle\hat{L},\alpha,\mu\rangle$, we see that the map $\alpha\co \hat{L}A\to\hat{L}\hat{L}A$ is a retract of $\alpha\co \hat{L}\hat{L}A\to\hat{L}\hat{L}\hat{L}A$.  Hence, if the first map is not an equivalence, then the second is not.
\end{proof}

\subsection{The discrete case}
\label{sec:8.7}
Let $\M$ be a bicomplete pointed category, viewed as a discrete model category (\ref{sec:4.3}), with a class $\G$ of injective models in $\Ho\M=\M$, and let $\I\subset\M$ be the full subcategory of $\G$--injective objects in $\M$.  By Lemma \ref{lem:8.8} below for $A\in\M$, there is a natural isomorphism
$$
\hat{L}_{\G}A~\iso~\lim_{f\co A\to I}I\in\M
$$
where $f$ ranges over the comma category $A\downarrow\I$, and the $\G$--completion $\alpha\co A\to\hat{L}_{\G}A$ is the canonical map to this limit.  Hence, $\hat{L}_{\G}\co \M\to\M$ is a right Kan extension of the inclusion functor $\I\to\M$ along itself, and may therefore be viewed as a \emph{codensity triple} functor (see \cite[X.7]{Mac}).  We have used the following:

\begin{lem}
\label{lem:8.8}
For $A\in\M$, there is a natural isomorphism $\hat{L}_{\G}A\iso\lim_{f\co A\to I}I$ where $f$ ranges over $A\downarrow\I$.
\end{lem}

\begin{proof}
Let $\alpha\co A\to J^\bullet$ be a $\G$--resolution of $A$ in $\M$.  Then $J^n\in\I$ for $n\geq 0$ by Lemma \ref{lem:5.3}, and 
$\alpha^*\co \Hom(J^\bullet,I)\to\Hom(A,I)$ is a trivial fibration in $\Scal$ for each $I\in\I$ by Corollary \ref{cor:3.16}.  Thus the maps $\alpha\co A\to J^0$ and $d^0,d^1\co J^0\to J^1$ satisfy the conditions: (i) $J^0,J^1\in\I$; (ii) $d^0\alpha=d^1\alpha$; (iii) if $f\co A\to I\in\I$, then there exists $\bar{f}\co J^0\to I$ with $\bar{f}\alpha=f$; and (iv) if $g_0,g_1\co J^0\to I\in\I$ and $g_0\alpha=g_1\alpha$, then there exists $\bar{g}\co J^1\to I$ with $\bar{g}d^0=g_0$ and $\bar{g}d^1=g_1$.  Hence, $\lim_{f\co A\to I}I$ is the equalizer of $d^0,d^1\co J^0\to J^1$, which is isomorphic to 
$\Tot J^\bullet\iso\hat{L}_{\G}A$.
\end{proof}

\subsection{The Bousfield--Kan case with an erratum} 
\label{sec:8.9}
By \ref{sec:7.7} and Corollary \ref{cor:8.2}, the Bousfield--Kan $R$--completion $\alpha\co X\to R_\infty X$ belongs to a triple on $\Ho_*$.  However, we no longer believe that it belongs to a triple  on $\Scal_*$ or $\Scal$, as we claimed in \cite[page 26]{BK}.  In that work, we correctly constructed functors $R_s\co \Scal_*\to\Scal_*$ with compatible transformations $1\to R_s$ and $R_sR_s\to R_s$ satisfying the left and right identity conditions for $0\leq s\leq\infty$, but we now think that our transformation $R_sR_s\to R_s$ is probably nonassociative for $s\geq 2$, because the underlying cosimplicial pairing $c$ in \cite[page 28]{BK} is nonassociative in cosimplicial dimensions $\geq 2$.  The difficulty arises because our ``twist maps''
do not compose to give actual summetric group actions on the $n$-fold composites $R\cdots R$ for $n\geq 3$.  The partial failure of our triple lemma in \cite{BK} does not seem to invalidate any of our other results, and new work of Libman \cite{Lib+} on homotopy limits for coaugmented functors shows that the functors $R_s$ must all still belong to triples on the homotopy category $\Ho_*$.  
%Moreover, A. Libman 
%\cite{Lib+} has recently shown that the Bendersky--Thompson tower functors with respect to an $S$--algebra spectrum all %belong to triples on $\Ho_*$, confirming that the functors $R_s$ belong to triples on $\Ho_*$ for $0\leq s\leq\infty$.

\section{Comparing different completions}
\label{sec:9}
We develop machinery for comparing different completion functors and apply it to show that the Bendersky--Thompson completions with respect to connective ring spectra are equivalent to Bousfield--Kan completions with respect to solid rings, although the associated homotopy spectral sequences may be quite different.   We continue to let $\C$ be a left proper, bicomplete, pointed simplicial model category with a class $\G$ of injective models in $\Ho\C$.   In addition, we suppose that $\C$ is factored and that $\G$ is functorial, so that the model category $\cC^{\G}$ is also factored by \ref{sec:4.2}.  Thus the $\G$--completion functor $\hat{L}_{\G}$ is defined on $\C$ (not just $\Ho\C$) by \ref{sec:5.7}.  We start by expressing the total derived functor $\Tot_{\G}\co \Ho(\cC^{\G})\to\Ho\C$ of \ref{cor:8.2} in terms of the prolonged functor $\hat{L}_{\G}\co \cC\to\cC$ with $(\hat{L}_{\G}X^\bullet)^n=\hat{L}_{\G}(X^n)$ for $n\geq 0$ 
and the homotopical $\Tot$ functor  $\underline{\Tot}\co \cC\to\C$ with $\underline{\Tot}X^\bullet=\Tot\underline X^\bullet$, where $\underline X^\bullet$ is a functorial Reedy fibrant replacement of $X^\bullet\in\cC$.

\begin{thm}
\label{thm:9.1}
For $Y^\bullet\in\cC$, there is a natural equivalence $$\Tot_{\G}Y^\bullet\homeq\underline{\Tot}(\hat{L}_{\G}Y^\bullet)$$ in $\Ho\C$.
\end{thm}

\begin{proof}
As in \ref{sec:6.11}, let \ $\con(Y^\bullet)\to\vec{Y}^{\bullet\bullet}$ be the functorial Reedy $\G$--resolution of \ $\con(Y^\bullet)$.  This induces a $\G$--equivalence of diagonals $Y^\bullet\to\Bar{\Bar{Y}}^\bullet$ with $\Bar{\Bar{Y}}^\bullet$ $\G$--fibrant and therefore induces 
$$
\Tot_{\G}Y^\bullet~\homeq~\Tot\Bar{\Bar{Y}}^\bullet~=~\Tot\diag\vec{Y}^{\bullet\bullet}~\homeq~\Tot\Tot^v\vec{Y}^{\bullet\bullet}~\homeq~\underline{\Tot}\Tot^v\vec{Y}^{\bullet\bullet}
$$
in $\Ho\C$ by Lemma \ref{lem:6.15}.  Now let $Y^n\to\bar{Y}^{n\bullet}$ be the functorial $\G$--resolution of $Y^n$ for $n\geq 0$, and functorially factor \ $\con(Y^\bullet)\to\bar{Y}^{\bullet\bullet}$ into a Reedy $\G$--trivial cofibration \ $\con(Y^\bullet)\to K^{\bullet\bullet}$ and a Reedy $\G$--fibration $K^{\bullet\bullet}\to\bar{Y}^{\bullet\bullet}$.
Next choose a map $K^{\bullet\bullet}\to\vec{Y}^{\bullet\bullet}$ extending \ $\con(Y^\bullet)\to\vec{Y}^{\bullet\bullet}$.  Since the maps \ $Y^n\to\vec{Y}^{n\bullet}$, \ $Y^n\to\bar{Y}^{n\bullet}$, \ and \ $Y^n\to K^{n\bullet}$ \ are $\G$--resolutions for $n\geq 0$, the maps $\vec{Y}^{n\bullet}\leftarrow K^{n\bullet}\to\bar{Y}^{n\bullet}$ are Tot-equivalences, and we obtain $\underline{\Tot}$--equivalences
$$
\Tot^v\vec{Y}^{\bullet\bullet} \xleftarrow{\qquad} \Tot^vK^{\bullet\bullet} \xrightarrow{\qquad}\Tot^v\bar{Y}^{\bullet\bullet}~\homeq~\hat{L}_{\G}Y^\bullet
$$
which combine to give $\underline{\Tot}\Tot^v\vec{Y}^{\bullet\bullet}\homeq\underline{\Tot}(\hat{L}_{\G}Y^\bullet)$ in $\Ho\C$.  This completes our chain of equivalences from \ $\Tot_{\G}Y^\bullet$ to \ $\underline{\Tot}(\hat{L}_{\G}Y^\bullet)$.
\end{proof}

\begin{cor}
\label{cor:9.2}
A $\G$--equivalence $X^\bullet\to Y^\bullet$ in $\cC$ induces an equivalence 
\ $\underline{\Tot}(\hat{L}_{\G}X^\bullet)\homeq\underline{\Tot}(\hat{L}_{\G}Y^\bullet)$ in $\Ho\C$.
\end{cor}

This follows immediately from Theorem \ref{thm:9.1} and specializes to give the following: 

\begin{cor}
\label{cor:9.3}
For an object $A\in\C$, each $\G$--equivalence $A\to Y^\bullet$ in  $\cC$ induces an equivalence $\hat{L}_{\G}A\homeq\underline{\Tot}(\hat{L}_{\G}Y^\bullet)$ in $\Ho\C$.
\end{cor}

\begin{defn}
\label{defn:9.4}
A \emph{$\G$--complete expansion} of an object $A\in\C$ consists of a $\G$--equivalence $A\to Y^\bullet$ in $\cC$ such that $Y^n$ is $\G$--complete for $n\geq 0$.
\end{defn}

Each weak $\G$--resolution of $A$ is a $\G$--complete expansion of $A$, and the completion part of Theorem \ref{thm:6.5} now generalizes to  the following:

\begin{thm}
\label{thm:9.5}
If $A\to Y^\bullet$ is a $\G$--complete expansion of an object $A\in\C$, then there is a natural equivalence $\hat{L}_{\G}A\homeq\underline{\Tot}Y^\bullet$ in $\Ho\C$.
\end{thm}

\begin{proof}
By Corollary \ref{cor:9.3}, the maps 
$\hat{L}_{\G}A\longrightarrow\underline{\Tot}\hat{L}_{\G}Y^\bullet\longleftarrow\underline{\Tot}Y^\bullet$
are weak equivalences in $\C$.
\end{proof}

%Old 9.6 omitted

By this theorem, any functorial $\G$--complete expansion of the objects in $\C$ gives a $\G$--completion functor on $\C$ which is ``essentially equivalent'' to $\hat{L}_{\G}$ since it is related to $\hat{L}_{\G}$ by natural weak equivalences.  The following theorem will show that different choices of $\G$ may give equivalent $\G$--completion functors even when they give very different $\G$--homotopy spectral sequences.

\begin{thm}
\label{thm:9.7}
Suppose $\G$ and $\G^\prime$ are classes of injective models in $\Ho\C$.  If each $\G$--injective object is $\G^\prime$--injective and each $\G^\prime$--injective object is $\G$--complete, then there is a natural equivalence $\hat{L}_{\G}A\homeq\hat{L}_{\G^\prime}A$ for $A\in\C$.
\end{thm}

\begin{proof}
Let $A\to J^\bullet$ be a $\G^\prime$--resolution of $A$.  Then $A\to J^\bullet$ is a $\G$--trivial cofibration by Corollary \ref{cor:3.16}, and $J^\bullet$ is termwise $\G^\prime$--injective.  Hence $A\to J^\bullet$ is a $\G$--complete expansion of $A$, and $\hat{L}_{\G}A\homeq\underline{\Tot}J^\bullet\homeq\hat{L}_{\G^\prime}A$ by Theorem \ref{thm:9.5}.
\end{proof}

For example, consider the Bendersky--Thompson completion $A\to\hat{A}_E$ of a space $A$ with respect to a ring spectrum $E$ as in \ref{sec:7.8}.  Suppose $E$ is \emph{connective} (ie, $\pi_i E=0$ for $i<0$), and suppose the ring $\pi_0E$ is commutative.  Let $R=\core(\pi_0E)$ be the subring 
$$
R~=~\{r\in\pi_0E~|~r\tensor1=1\tensor r\in\pi_0E\tensor\pi_0E\},
$$
and recall that $R$ is \emph{solid} (ie, the multiplication $R\tensor R\to R$ is an isomorphism) by \cite{BK1972}.

\begin{thm}
\label{thm:9.8}
If $E$ is a connective ring spectrum with commutative $\pi_0 E$, then there are natural equivalences $\hat{A}_E\homeq(\pi_0 E)_\infty A\homeq R_\infty A$ for $A\in\Ho_*$ where $R=\core(\pi_0 E)$.
\end{thm}

\begin{proof}
Let $\G^\prime$ (resp.\ $\G$) be the class of all $\Omega^\infty N\in\Ho_*$ for $E$--module (resp.\ $H\pi_0 E$--module) spectra $N$.  Then $\G\subset\G^\prime$ since each $H\pi_0 E$--module spectrum is an $E$--module spectrum via the map $E\to H\pi_0 E$, and hence each $\G$--injective space is $\G'$--injective.  If $N$ is an $E$--module spectrum, then $(\pi_0 E)_\infty\Omega^\infty N\homeq\Omega^\infty N$ by \cite[II.2]{BK}.  Hence each $\G^\prime$--injective space $J$ is $\G$--complete, since it is a retract of $\Omega^\infty N$ for $N=E\wedge\Sigma^\infty J$.  Consequently, $\hat{A}_E\homeq\hat{L}_{\G^\prime}A\homeq\hat{L}_{\G}A\homeq(\pi_0 E)_\infty A$ by Theorem \ref{thm:9.7}, and $(\pi_0 E)_\infty A\homeq R_\infty A$ by \cite[page 23]{BK}.
\end{proof}

\subsection{Examples of $E$--completions}
\label{sec:9.9}
In \cite{BK1972} and \cite[6.4]{Bou1979}, we determined all solid rings $R$, and they are: (I) $R=\Z[J^{-1}]$ for a set $J$ of primes; (II) $R=\Z/n$ for $n\geq2$; (III) $R=\Z[J^{-1}]\times\Z/n$ for $n\geq2$ and a set $J$ of primes including the factors of $n$; and (IV) $R=\core(\Z[J^{-1}]\times\prod_{p\in K}\Z/p^{e(p)})$ for infinite sets $K\subset J$ of primes and positive integers $e(p)$.  In \cite[I.9]{BK}, we showed that the completions $R_\infty X$ in cases (I)--(III) can be expressed as products of their constituent completions $\Z[J^{-1}]_\infty X$ and $(\Z/p)_\infty X$ for the prime factors $p$ of $n$, and we extensively studied these basic completions.  We found that a nilpotent space $X$ is always $R$--good in cases (I) and (II), but is ``usually'' $R$--bad in cases (III) and (IV).  For instance $K(\Z,m)$ for $m\geq1$ is  $R$--bad in cases (III) and (IV).  These results are now applicable to the completion $\hat{X}_E$ of a space $X\in\Ho_*$ with respect to a connective ring spectrum $E$ with $\pi_0E$ commutative.  For instance, we have $\hat{X}_E\homeq\Z_\infty X$ for $E=S$ and $\hat{X}_E\homeq\Z_{(p)\infty}X$ for $E=BP$.

\subsection{The loop-suspension completions}
\label{sec:9.10}
We may also apply Theorem \ref{thm:9.7} to reprove the result that loop-suspension completions of spaces are equivalent to $\Z$--completions.  In more detail, for a fixed integer $n\geq 1$, we consider the $n$-th loop-suspension completion (\ref{sec:7.9}) of a space $A\in\Scal_*$ given by $\hat{L}_{\G}A$ where $\G=\{\Omega^nY~|~Y\in\Ho_*\}$, and we compare it with the Bousfield--Kan $\Z$--completion $\hat{L}_{\Hcal}A\homeq\Z_\infty A$ where $\Hcal=\{\Omega^\infty N~|~N \text{ is an $H$--module spectrum}\}$.  Since the $\G$--injective spaces are the retracts of the $n$-fold loop spaces, they have nilpotent components and are $\Hcal$--complete.  Thus, since $\Hcal\subset \G$, Theorem \ref{thm:9.7} shows 
$\hat{L}_{\G}A\homeq\hat{L}_{\Hcal}A$, and the $n$th loop-suspension completion of $A$ is equivalent to $\Z_\infty A$.

\section{Bendersky-Thompson completions of fiber squares}
\label{sec:10}
Let $\C$ remain a left proper, bicomplete, pointed simplicial model category with a class $\G$ of injective models in $\Ho\C$.  Also suppose that $\C$ is factored and $\G$ is functorial so that the $\G$--completion functor $\hat{L}_{\G}$ is defined on $\C$ (not just $\Ho\C$) by \ref{sec:5.7}.  In this section, we show that $\hat{L}_{\G}$ preserves fiber squares whose ``$\G$--cohomology cobar spectral sequences collapse strongly,''  and we specialize this result to the Bendersky--Thompson completions (see Theorems \ref{thm:10.11} and \ref{thm:10.12}).   We need a weak assumpton on the following: 

\subsection{Smash products in $\Ho\C$}
\label{sec:10.1}
For $A,B\in\Ho\C$, let $A\wedge B\in\Ho\C$ be the smash product represented by the homotopy cofiber of the coproduct-to-product map $A\vee B\to A\times B$ for cofibrant-fibrant objects $A,B\in\C$.  We assume that the functor 
$-\wedge B:\Ho\C\to\Ho\C$ has a right adjoint $(-)^B:\Ho\C\to\Ho\C$.  This holds as usual in $\Ho\Scal_*=\Ho_*$, and it is easy to show the following:

\begin{lem}
\label{lem:10.2}
For an object $B\in\Ho\C$, the following are equivalent:
\begin{enumerate}
\item[\rm(i)]if a map $X\to Y$ in $\Ho\C$ is $\G$--monic, then so is $X\wedge B\to Y\wedge B$;
\item[\rm(ii)]if an object $I\in\Ho\C$ is $\G$--injective, then so is $I^B$;
\item[\rm(iii)]for each $G\in\G$ and $i\geq 0$, the object $(\Omega^iG)^B$ is $\G$--injective.
\end{enumerate}
\end{lem}

\begin{defn}
\label{defn:10.3}
An object $B\in\Ho\C$ will be called \emph{$\G$--flat} (for smash products) when it satisfies the equivalent conditions of Lemma \ref{lem:10.2}.  An object $B\in\C$ (resp.\ $B^\bullet\in\cC$) will also be called \emph{$\G$--flat} when $B$ (resp.\ each $B^n$) is $\G$--flat in $\Ho\C$.
\end{defn}

\begin{lem}
\label{lem:10.4}
If $f\co X^\bullet\to Y^\bullet$ and $g\co B^\bullet\to C^\bullet$ are $\G$--equivalences of termwise fibrant objects in $\cC$ such that $Y^\bullet$ and $B^\bullet$ are $\G$--flat, then  $f\times g\co X^\bullet\times B^\bullet\to Y^\bullet\times C^\bullet$ is also a $\G$--equivalence.
\end{lem}

\begin{proof}
Working in $c(\Ho\C)$ instead of $\cC$, we note that $f\wedge B^n\co X^\bullet\wedge B^n\to Y^\bullet\wedge B^n$ is a $\G$--equivalence for $n\geq 0$ by Lemma \ref{lem:10.5} below, since $(\Omega^iG)^{B^n}\in\Ho\C$ is a $\G$--injective group object for each $G\in\G$ and $i\geq 0$.  Hence, $f\wedge B^\bullet\co X^\bullet\wedge B^\bullet\to Y^\bullet\wedge B^\bullet$ is a $\G$--equivalence as in the proof of Lemma \ref{lem:6.9}.  Similarly $Y^\bullet\wedge g\co Y^\bullet\wedge B^\bullet\to Y^\bullet\wedge C^\bullet$ is a $\G$--equivalence, and hence so is $f\wedge g\co X^\bullet\wedge B^\bullet\to Y^\bullet\wedge C^\bullet$. Thus the ladder
$$
\begin{CD}
X^\bullet\vee B^\bullet  @>>>  X^\bullet\times B^\bullet  @>>>  X^\bullet\wedge B^\bullet  \\
@VV{f\vee g}V                   @VV{f\times g}V                   @VV{f\wedge g}V     \\
Y^\bullet\vee C^\bullet  @>>>  Y^\bullet\times C^\bullet  @>>>  Y^\bullet\wedge C^\bullet
\end{CD}
$$
is carried by $[-,G]_*$ to a ladder of short exact sequences of simplicial groups such that 
$(f\vee g)^*$ and $(f\wedge g)^*$ are weak equivalences.   Consequently $(f\times g)^*$ is a weak equivalence.
\end{proof}

We have used the following:

\begin{lem}
\label{lem:10.5}
If $f\co X^\bullet\to Y^\bullet\in\cC$ is a $\G$--equivalence and $I\in\Ho\C$ is a $\G$--injective group object, then $f^*\co [Y^\bullet,I]_*\to[X^\bullet,I]_*$ is a weak equivalence of simplicial groups.
\end{lem}

\begin{proof}
The  class of $\G$--monic maps in $\Ho\C$ is clearly the same as the class of $\G^\prime$--monic maps for $\G^\prime=\G\cup\{I\}$.  Hence, $\G$ and $\G^\prime$ give the same model category structure on $\cC$ by \ref{sec:4.1}, and $f\co X^\bullet\to Y^\bullet$ is a $\G^\prime$--equivalence in $\cC$.
\end {proof}

\begin{thm}
\label{thm:10.6}
Suppose the $\G$--injectives in $\Ho\C$ are $\G$--flat.  If $A,B,M\in\Ho\C$ are objects with $A$ or $B$ $\G$--flat, then there is  a natural equivalence $\hat{L}_{\G}(A\times B)\homeq\hat{L}_{\G}A\times\hat{L}_{\G}B$ and a natural isomorphism
$$
E^{s,t}_r(A\times B;M)_{\G}~\homeq~E^{s,t}_r(A;M)_{\G}\times E^{s,t}_r(B;M)_{\G}
$$
for $2\leq r\leq\infty +$ and $0\leq s\leq t$.
\end{thm}

\begin{proof}
We may suppose  $A$ and $B$ are fibrant in $\C$ and take $\G$--resolutions $A\to\bar{A}^\bullet$ and $B\to\bar{B}^\bullet$ in $\cC$.  Then the product $A\times B\to\bar{A}^\bullet\times\bar{B}^\bullet$ is  a weak $\G$--resolution by Lemma \ref{lem:10.4}, and the result follows from Theorem \ref{thm:6.5}.
\end{proof}

We now study the action of $\hat{L}_{\G}$ on a commutative square
\begin{equation}
\label{eq:10.7}
\begin{CD}
C    @>>>    B  \\
@VVV        @VVV\\
A    @>>>  \Lambda
\end{CD}
\end{equation}
of fibrant objects in $\C$ using the following:

\subsection{The geometric cobar construction}
\label{sec:10.8}
Let $\B(A,\Lambda,B)^\bullet\in\cC$ be the usual geometric cobar construction with
$$
\B(A,\Lambda,B)^n~=~A\times\Lambda\times\cdots\times\Lambda\times B
$$
for $n\geq 0$ where the factor $\Lambda$ occurs $n$ times (see \cite{Rec}).  It is straightforward to show that $\B(A,\Lambda,B)^\bullet$ is Reedy fibrant with
$$
\Tot\B(A,\Lambda,B)^\bullet~\iso~ P(A,\Lambda,B)
$$
where $P(A,\Lambda,B)$ is the double mapping path object defined by the pullback
$$
\begin{CD}
P(A,\Lambda,B)   @>>>   \hom(\Delta^1,\Lambda)\\
@VVV                      @VVV \\
A\times B        @>>>   \Lambda\times\Lambda.
\end{CD}
$$
Thus $P(A,\Lambda,B)$ represents the homotopy pullback of the diagram $A\to\Lambda\leftarrow B$ (see \cite[\S10]{DS}),
 and \eqref{eq:10.7} is called a \emph{homotopy fiber square} when the map $C\to P(A,\Lambda,B)$ is a weak equivalence.

Our main fiber square theorem for $\G$--completions is the following:

\begin{thm}
\label{thm:10.9}
Suppose the $\G$--injectives in $\Ho\C$ are $\G$--flat.  If \eqref{eq:10.7} is a  square of $\G$--flat fibrant objects such that the augmentation $C\to\B(A,\Lambda,B)^\bullet$ is a $\G$--equivalence, then $\hat{L}_{\G}$ carries \eqref{eq:10.7} to a homotopy fiber square.
\end{thm}

\begin{proof}
Since $C\to\B(A,\Lambda,B)^\bullet$ \ is a $\G$--equivalence, it induces an equivalence $\hat{L}_{\G}C\homeq\underline{\Tot}\hat{L}_{\G}\B(A,\Lambda,B)^\bullet$ by Corollary \ref{cor:9.3}, and there are equivalences 
$$
\underline{\Tot}\hat{L}_{\G}\B(A,\Lambda,B)^\bullet
~\homeq~\underline{\Tot}\B(\hat{L}_{\G}A,\hat{L}_{\G}\Lambda,\hat{L}_{\G}B)^\bullet
~\homeq~P(\hat{L}_{\G}A,\hat{L}_{\G}\Lambda,\hat{L}_{\G}B)
$$
by Theorem \ref{thm:10.6}.  Hence, $\hat{L}_{\G}C$ is equivalent to the homotopy pullback of 
$\hat{L}_{\G}A\to\hat{L}_{\G}\Lambda\leftarrow\hat{L}_{\G}B$.
\end{proof}

The hypothesis that the augmentation $C\to\B(A,\Lambda,B)^\bullet$ is a $\G$--equivalence may be reformulated to say that \emph{the $\G$--cohomology cobar spectral sequences collapse strongly} for \eqref{eq:10.7}, although we shall not develop that viewpoint here.

\subsection{The Bendersky--Thompson case}
\label{sec:10.10}
For a commutative ring spectrum $E$, we consider the \emph{Bendersky--Thompson $E$--completion} $A\to\hat{A}_E=\hat{L}_{\G}A$ of a space $A\in\Scal_*$ with respect to the class of injective models
$$
\G~=~\{\Omega^\infty N~|~N\text{ is a $E$--module spectrum}\}~\subset~\Ho_*
$$
as in \ref{sec:7.8}.  All spaces in $\Ho_*$ are now $\G$--flat, and Theorem \ref{thm:10.9} will apply to the square \eqref{eq:10.7} provided that \emph{the $N^*$--cobar spectral sequence collapses strongly} for each $E$--module spectrum $N$ in the sense that 
$$
\pi_sN^*\B(A,\Lambda,B)^\bullet ~\iso~
\begin{cases}
N^*C   &\text{for $s=0$}\\
0      &\text{for $s>0$.}
\end{cases}
$$
Here we may assume that $N$ is an extended $E$--module spectrum since any $N$ is a homotopy retract of $E\wedge N$.  To eliminate $N$ from our hypotheses, we suppose:
\begin{enumerate}
\item[(i)]$E$ satisfies the \emph{Adams UCT condition} namely that the map $N^*X \to$\break $\Hom^*_{E_*}(E_*X,\pi_*N)$ is an isomorphism for  each $X\in\Ho_*$ with $E_*X$ projective over $E_*$ and each extended $E$--module spectrum $N$;  
\item[(ii)]$E_*A$, $E_*\Lambda$, $E_*B$ and $E_*C$ are projective over $E_*$.
\end{enumerate}
Condition (i) holds for many common ring spectra $E$, including  the $p$--local ring spectrum $K$ and arbitrary $S$--algebras by \cite[page 284]{Ada}  and \cite[page 82]{EKMM}.  Condition (ii) implies that 
$$
E_*\B(A,\Lambda,B)^n~\iso~ E_*A\tensor_{E_*}E_*\Lambda\tensor_{E_*}\cdots\tensor_{E_*}E_*B
$$
is projective over $E_*$ for $n\geq 0$, and we say that \emph{the $E_*$--cobar spectral sequence collapses strongly} when 
$E_*C\to E_*\B(A,\Lambda,B)^\bullet$ is split exact as a complex over $E_*$.  Now Theorem \ref{thm:10.9}  implies the following:

\begin{thm}
\label{thm:10.11}
Suppose $E$ is a commutative ring spectrum satisfying the Adams UCT-condition.  If the spaces of \eqref{eq:10.7} have $E_*$--projective homologies and the $E_*$--cobar spectral sequence collapses strongly, then the Bendersky--Thompson $E$--completion  functor carries \eqref{eq:10.7} to a homotopy fiber square.
\end{thm}

Specializing this to $E=K$, we suppose that the spaces of \eqref{eq:10.7} have $K_*$--free homologies, and we say that \emph{the $K_*$--cobar spectral sequence collapses strongly} if
$$
\Cotor^{K_*\Lambda}_s(K_*A,K_*B)~=~
\begin{cases}
K_*C        &\text{for $s=0$}\\
0           &\text{for $s>0$.}
\end{cases}
$$
Now Theorem \ref{thm:10.11} reduces to the following:

\begin{thm}
\label{thm:10.12}
If the spaces of \eqref{eq:10.7} have $K_*$--free homologies and the $K_*$--cobar specral sequence collapses strongly, then the Bendersky--Thompson $K$--completion functor carries \eqref{eq:10.7} to a homotopy fiber square.
\end{thm}

This result is applied by Bendersky and Davis in \cite{BD2001}.

\section{$p$--adic $K$--completions of fiber squares}
\label{sec:11}

Working at an arbitrary prime $p$, we now consider a $p$--adic variant of the Bendersky--Thompson $K$--completion of spaces and establish an improved fiber square theorem for it.  We also briefly consider the associated homotopy spectral sequence which seems especially applicable to spaces whose $p$--adic $K$--cohomologies are torsion-free with Steenrod-Epstein-like $U(M)$ structures as in \cite{Bou1996a}.  We first recall the following:

\subsection{The $p$--completion of a space or spectrum}
\label{sec:11.1}
For a space $A\in\Scal_*$, we let $\hat{A}=A_{H/p}$ be the \emph{$p$--completion} given by the $H/p_*$--localization of 
\cite{Bou1975}.  This is equivalent to the $S/p_*$--localization and, when $A$ is nilpotent, is equivalent to the $p$--completion $(Z/p)_\infty A$ of \cite{BK}.  For a spectrum $E$, we likewise let $\hat{E}=E_{S/p}$ be the \emph{$p$--completion} given by the $S/p_*$--localization of \cite{Bou1979}.  Thus, when the groups $\pi_*E$ are finitely generated, we have $\pi_*\hat{E}=\pi_*E\tensor\hat{\Z}_p$ using the $p$--adic integers  $\hat{\Z}_p$.  We now introduce the following:

\subsection{The $p$--adic $K$--completion}
\label{sec:11.2}
The triple on $\Ho_*$ carrying a space $X$ to $\Omega^\infty(K\wedge\Sigma^\infty X)\widehat{~}$ satisfies the conditions of  \ref{sec:7.5} and thus determines a class of injective models
$$
\hat{\G}~=~\{\Omega^\infty(K\wedge\Sigma^\infty X)\widehat{~}~|~X\in\Ho_*\}~\subset~\Ho_*.
$$
For spaces $A,M\in\Scal_*$, let $\hat{A}_{\hat{K}}=\hat{L}_{\hat{\G}}A$ be the resulting \emph{$p$--adic $K$--completion} and consider the associated homotopy spectral sequence $\{E^{s,t}_r(A;M)_{\hat{K}}\}=\{E^{s,t}_r(A;M)_{\hat{\G}}\}$.  We could equivalently use the class of injective models
$$
\hat{\G}'~=~\{\Omega^\infty N~|~N\text{ is a $p$--complete $K$--module spectrum}\}~\subset~\Ho_*
$$
or less obviously, when $K^*(A;\hat{\Z}_p)$ is torsion-free, use the class of injective models representing the $p$--adic $K$--cohomology theory $K^*(-;\hat{\Z}_p)$ as in \ref{sec:4.6}.

\subsection{Comparison with the Bendersky--Thompson $K$--completion}
\label{sec:11.3}
For the $p$--local ring spectrum $K$ and a space $A\in\Scal_*$, let $\hat{A}_K=\hat{L}_{\G}A$ be the Bendersky--Thompson $K$--completion obtained using the class of injective models 
$$
\G~=~\{\Omega^\infty N~|~N\text{ is a $K$--module spectrum}\}~\subset~\Ho_*
$$
as in \ref{sec:7.8} or \ref{sec:10.10}.  Also consider the associated homotopy spectral sequence  $\{E^{s,t}_r(A;M)_K\}=\{E^{s,t}_r(A;M)_{\G}\}$ for $A,M\in\Scal_*$.
Since $\hat{\G}\subset\G$, there is a natural map $\hat{A}_K\to\hat{A}_{\hat{K}}$ constructed as follows for a space $A\in\Scal_*$.  First take a $\G$--resolution $A\to I^\bullet$ of $A$ and then take a $\hat{\G}$--resolution $I^\bullet\to J^\bullet$ of $I^\bullet$ in $\cS_*$.  Since the composed map $A\to J^\bullet$ is a $\hat{\G}$--resolution of $A$, the map $I^\bullet\to J^\bullet$ induces the desired map $\hat{A}_K\homeq\Tot I^\bullet\to\Tot J^\bullet\homeq\hat{A}_{\hat{K}}$.  It also induces a map $\{E^{s,t}_r(A;M)_K\} \to \{E^{s,t}_r(A;M)_{\hat{K}}\}$ of homotopy spectral sequences for $A,M\in\Scal_*$.  The following theorem will show that these maps are ``almost $p$--adic equivalences.''   For a space $Y\in\Ho_*$, let 
$Y\langle n\rangle\in\Ho_*$ be the $(n-1)$--connected section of $Y$, and let $Y\langle\tilde{n}\rangle\in\Ho_*$ be the section with
$$
\pi_i Y\langle\tilde{n}\rangle~=~
\begin{cases} 
\pi_i Y          &\text{for $i>n$}\\
(\pi_n Y)\tilde{\ }  &\text{for $i=n$}\\
0                &\text{for $i<n$}
\end{cases}
$$
where $(\pi_n Y)\tilde{\ }$ is the divisible part of $\pi_n Y$ assuming $n\geq 2$.

\begin{thm}
\label{thm:11.4}
If $A,M\in\Scal_*$ are spaces with $\tilde{H}_*(M;Q)=0$, then:
\begin{enumerate}
\item[\rm(i)]$\hat{A}_{\hat{K}}\langle 3\rangle$ is the $p$--completion of $\hat{A}_K\langle\tilde{2}\rangle$;
\item[\rm(ii)]$[M,\hat{A}_K]_*\iso[M,\hat{A}_{\hat{K}}]_*$;
\item[\rm(iii)]$E^{s,t}_r(A;M)_K\iso E^{s,t}_r(A;M)_{\hat{K}}$ for $0\leq s\leq t$ and $2\leq r\leq\infty+$.
\end{enumerate}
\end{thm}

This will be proved in \ref{sec:11.10}.  For a space $A$, we may actually construct the $p$--adic $K$--completion of $A$ and the associated homotopy spectral sequence quite directly from the Bendersky--Thompson triple resolution $A\to K^\bullet A$ of \ref{sec:7.8}.  We simply apply the $p$--completion functor to give a map $A\to\widehat{K^\bullet A}$ in $\cS_*$ and obtain the following: 

\begin{thm}
\label{thm:11.5}
For a space $A\in\Scal_*$, the map $A\to\widehat{K^\bullet A}$ is a weak $\hat{\G}$--resolution of $A$.  Hence $\hat{A}_{\hat{K}}\homeq\underline{\Tot}(\widehat{K^\bullet A})$  and
\ $E^{s,t}_r(A;M)_{\hat{K}}\iso E^{s,t}_r(\widehat{K^\bullet A};M)$\  for $M\in\Scal_*$, $0\leq s\leq t$, and 
$2\leq r\leq \infty+$.
\end{thm}

This will be proved in \ref{sec:11.9}.  We now turn to our fiber square theorem for the $p$--adic $K$--completion.  For a commutative square of fibrant spaces
\begin{equation}
\label{eq:11.6}
\begin{CD}
C     @>>>     B\\
@VVV          @VVV\\
A     @>>>    \Lambda
\end{CD}
\end{equation}
we say that the $K_*(-;Z/p)$--cobar spectral sequence \emph{collapses strongly} when 
$$
\Cotor^{K_*(\Lambda;\Z/p)}_s(K_*(A;\Z/p),K_*(B;\Z/p))~=~
\begin{cases}
K_*(C;\Z/p)   &\text{for $s=0$}\\
0             &\text{otherwise.}
\end{cases}
$$

\begin{thm}
\label{thm:11.7}
If the spaces in \eqref{eq:11.6} have torsion-free $K^*(-;\hat{\Z}_p)$--cohomol\-og\-ies and the $K_*(-;\Z/p)$--cobar spectral sequence collapses strongly, then the $p$--adic $K$--completion functor carries \eqref{eq:11.6} to a homotopy fiber square.
\end{thm}

This will be proved below in \ref{sec:11.12} using our general fiber square theorem (\ref{thm:10.9}).  It applies to a broader range of examples than its predecessor Theorem \ref{thm:10.12} for the Bendersky--Thompson $K$--completion, and we remark that its strong collapsing hypothesis holds automatically by \cite[Theorem 10.11]{Bou1996} whenever the spaces are connected and the coalgebra map $K_*(B;\Z/p)\to K_*(\Lambda;\Z/p)$ belongs to an epimorphism of graded bicommutative Hopf algebras (with possibly artificial multiplications).  We devote the rest of this section to proving the above theorems.

\begin{lem}
\label{lem:11.8}
If $N\in\Ho^s$ is a $K$--module spectrum, then the space $\widehat{\Omega^\infty N}$ is $\hat{\G}$--injective.
\end{lem}

\begin{proof}
The spaces $\widehat{\Omega^\infty N}\langle 1\rangle$ and $\Omega^\infty\hat{N}\langle 1\rangle$ can be expressed as 
$$
\widehat{\Omega^\infty N}\langle 1\rangle ~=~ SUJ_1\times UJ_2\times BUJ_3
$$
$$
\Omega^\infty\hat{N}\langle 1\rangle ~=~ UJ_1\times UJ_2\times BUJ_3
$$
for $\Ext$-$p$--complete abelian groups $J_1$,$J_2$,$J_3$ with $J_1=\Hom(\Z_{p^\infty},\pi_0 N)$ tors\-ion-free.  Since $SUJ_1$ is a retract of $UJ_1$ by \cite[Lemma 2.1]{Mis}, $\widehat{\Omega^\infty N}\langle 1\rangle$ is a retract of $\Omega^\infty\hat{N}\langle 1\rangle$, and both spaces are $\hat{\G}$--injective.  The lemma now follows since $\widehat{\Omega^\infty N}\homeq\widehat{\Omega^\infty N}\langle 1\rangle\times K(\pi_0N,0)$ and since   $K(\pi_0N,0)$ is also $\hat{\G}$--injective because it is discrete.
\end{proof}

\subsection{Proof of Theorem \ref{thm:11.5}}
\label{sec:11.9}
Since $A\to K^\bullet A$ is a $\G$--equivalence, it is also a $\hat{\G}$--equivalence, and hence so is 
$A\to\widehat{K^\bullet A}$.  Since the terms of $\widehat{K^\bullet A}$ are $\hat{\G}$--injective by Lemma \ref{lem:11.8}, this implies that  $A\to\widehat{K^\bullet A}$ is a weak $\hat{\G}$--resolution.  The final statement follows from Theorem \ref{thm:6.5}.
\endproof

\subsection{Proof of Theorem 11.4}
\label{sec:11.10}
For $0\leq s\leq\infty$, we obtain a homotopy fiber square
$$
\begin{CD}
\underline{\Tot}_s(K^\bullet A)    @>>>     \underline{\Tot}_s(\widehat{K^\bullet A})\\
@VVV                                            @VVV\\
\underline{\Tot}_s(K^\bullet A)_{(0)}   @>>>    \underline{\Tot}_s(\widehat{K^\bullet A})_{(0)}
\end{CD}
$$
by applying $\underline{\Tot}_s$ to the termwise arithmetic square \cite{DDK} of $K^\bullet A$.  Since the lower spaces of the square are $HQ_*$--local  \cite[page 192]{Bou1975}, the upper map has an $HQ_*$--local homotopy fiber and induces an equivalence 
$$
\map_*(M,\underline{\Tot}_s(K^\bullet A))~\homeq~\map_*(M,\underline{\Tot}_s(\widehat{K^\bullet A}))
$$
Thus by Theorem \ref{thm:11.5}, the map  $\hat{A}_K\to\hat{A}_{\hat{K}}$ has an $HQ_*$--local homotopy fiber and induces an equivalence $\map_*(M,\hat{A}_K)\homeq\map_*(M,\hat{A}_{\hat{K}})$.  The theorem now follows easily.  
\endproof

\begin{lem}
\label{lem:11.11}
For a space or spectrum $X$ with $K^*(X;\hat{\Z}_p)$ torsion-free and for an $\Ext$-$p$--complete abelian group $J$, the Pontrjagin dual $K^*(X;\hat{\Z}_p)^{\#}$ is divisible $p$--torsion with natural isomorphisms
$$
K_*(X;\Z/p)~\iso~K^*(X;\hat{\Z}_p)^{\#}\backslash p
$$
$$
K^*(X;J)~\iso~\Ext(K^*(X;\hat{\Z}_p)^{\#},J).
$$ 
\end{lem}

\begin{proof}
We can assume that $X$ is a spectrum and obtain natural isomorphisms
$$
K^*(X;\hat{\Z}_p)^{\#}~\iso~K_*(X;{\Z}/p^\infty)~\iso~K_{*-1}\tau_p X
$$
by \cite[Proposition 10.1]{Bou1999} where $\tau_p X$ is the $p$--torsion part of $X$.  Since these groups are divisible $p$--torsion and since $J$ is $\Ext$-$p$--complete, there are natural isomorphisms
$$
K^*(X;J)~\iso~K^*(\tau_p X;J)~\iso~\Ext(K_{*-1}\tau_p X,J)
$$
because $\Hom(K_*\tau_p X,J)=0$, and the lemma follows easily. 
\end{proof}

\subsection{Proof of Theorem \ref{thm:11.7}}
\label{sec:11.12}
Since all spaces in $\Ho_*$ are $\hat{\G}$--flat, it suffices by Theorem \ref{thm:10.9} to show that $C\to\B(A,\Lambda,B)^\bullet$ is a $\hat{\G}$--equivalence.  Since the augmented cochain complex 
$K_*(C;\Z/p)\to K_*(\B(A,\Lambda,B)^\bullet;\Z/p)$ is acyclic, the complex $K^*(C;\hat{\Z}_p)^{\#}\to K^*(\B(A,\Lambda,B)^\bullet;\hat{\Z}_p)^{\#}$ of divisible $p$--torsion groups must also be acyclic by Lemma \ref{lem:11.11}.  Hence, this complex must be contractible, and  the complex  $K^*(\B(A,\Lambda,B)^\bullet;J)\to K^*(C;J)$ must be acyclic for each $\Ext$-$p$--complete abel\-ian group $J$ by Lemma \ref{lem:11.11}.  Thus $C\to\B(A,\Lambda,B)^\bullet$ is a $\hat{\G}$--equivalence. \endproof

\section{The unpointed theory}
\label{sec:12}

As in \cite{GH}, much of the preceding work can be generalized to unpointed model categories.  In this section, we develop such a generalization (\ref{thm:12.4}) of the existence theorem (\ref{thm:3.3}) for $\G$--resolution model categories, and then briefly discuss the resulting unpointed theory of $\G$--resolutions, right derived functors, and $\G$--completions.  This leads, for instance, to unpointed Bendersky--Thompson completions of spaces.  We start with preliminaries on loop objects in unpointed model categories.

Let $\C$ be a model category with terminal object $e$, and let $\C_* = e\downarrow\C$ denote the associated pointed model category whose \emph{weak equivalences}, \emph{cofibrations}, and \emph{fibrations} are the maps having these properties when basepoints are forgotten.  The forgetful functor $\C_*\to\C$ is a Quillen right adjoint of the functor $\C\to\C_*$ sending $X\mapsto X\coprod e$ and has a total right derived functor $\Ho\C_*\to\Ho\C$ (see \ref{sec:4.7}).  We let 
$J\co \Ho\C_*\to(\Ho\C)_*$  be the associated functor to the pointed category $(\Ho\C)_* = [e]\downarrow\Ho\C$.

\begin{lem}
\label{lem:12.1}
For a  left proper model category $\C$, the isomorphism classes of objects in $\Ho\C_*$ correspond to the isomorphism classes of objects in $(\Ho\C)_*$ via the functor $J$.
\end{lem}

\begin{proof}
We first choose a trivial fibration $\check{e}\to e$ in $\C$ with $\check{e}$ cofibrant.  Then an object $X\in(\Ho\C)_*$ is represented by a cofibration $\check{e}\to X$ in $\C$, and the map $X\to X/\check{e}$ is a weak equivalence since $\C$ is left proper.  Hence $X\homeq J(X/\check{e})$ in $(\Ho\C)_*$.  For objects $W_1,W_2\in\Ho\C_*$ with $J(W_1)\homeq J(W_2)$, we may choose fibrant representatives $W_1,W_2\in\C_*$ and factor each $\check{e}\to e\to W_i$ into a cofibration $\check{e}\to\check{W}_i$ and a trivial fibration $\check{W}_i\to W_i$ in $\C$.  Using the homotopy extension theorem \cite[I.1.7]{Qui}  and the equivalence $J(W_1)\homeq J(W_2)$, we obtain a weak equivalence $\check{W}_1\to\check{W}_2$ under $\check{e}$.  Hence $W_1\homeq\check{W}_1/\check{e}\homeq \check{W}_2/\check{e}\homeq W_2$ in $\Ho\C_*$.
\end{proof}

\subsection{Loop objects in $(\Ho\C)_*$}
\label {sec:12.2}
For a left proper model category $\C$ and $n\geq 0$, the ordinary $n$--fold loop functor $\Omega^n\co \Ho\C_*\to\Ho\C_*$ now determines an $n$--fold loop operation $\Omega^n$ on the isomorphism classes of objects in $(\Ho\C)_*$ via the correspondence of Lemma \ref{lem:12.1}.  Thus for each object $Y\in(\Ho\C)_*$, we obtain an object $\Omega^nY\in(\Ho\C)_*$ defined up to isomorphism,  where  $\Omega^0Y=Y$.  We note that $\Omega^nY$ admits a group object structure in $\Ho\C$ for $n\geq 1$, which is  abelian for $n\geq 2$, since it comes from an $n$--fold loop object of $\Ho\C_*$ via a right adjoint functor $\Ho\C_*\to\Ho\C$.  For $X\in\Ho\C$, we let 
$$
[X,Y]_n~\iso~[X,\Omega^nY]~\iso~\Hom_{\Ho\C}(X,\Omega^nY)
$$
be the resulting homotopy set for $n=0$, group for $n=1$, or abelian group for $n\geq 2$.  When the original category $\C$ is pointed, we can identify $\C_*$ with $\C$, and our constructions give the usual objects $\Omega^nY\in\Ho\C$ and sets or groups $[X,Y]_n$.

\subsection{The $\G$--resolution model category}
\label{sec:12.3}
For a left proper model category $\C$, let $\G$ be a class of group objects in $\Ho\C$.  Then each $G\in\G$, with its unit map, represents an object of $(\Ho\C)_*$ and thus has an $n$--fold loop object $\Omega^nG\in\Ho\C$ giving an associated homotopy functor $[-,G]_n$ on $\Ho\C$ for $n\geq 0$.  A map $i\co A\to B$ in $\Ho\C$ is called 
\emph{$\G$--monic} when $i^*\co [B,G]_n\to[A,G]_n$ is onto for each $G\in\G$ and $n\geq 0$, and an object $Y\in\Ho\C$ is called \emph{$\G$--injective} when $i^*\co [B,Y]\to[A,Y]$ is onto for each $\G$--monic map $i\co A\to B$ in $\Ho\C$.  We retain the other definitions in \ref{sec:3.1} and \ref{sec:3.2}, and we obtain a structured simplicial category $\cC^{\G}$.  This leads to our most general existence theorem for resolution model categories.

\begin{thm}[after Dwyer--Kan--Stover]
\label{thm:12.4}
If $\C$ is a left proper model category with a class $\G$ of injective models in $\Ho\C$, then $\cC^{\G}$ is a left proper simplicial model category.
\end{thm}  

The proof proceeds exactly as in \ref{prop:3.4}--\ref{lem:3.23}, but thereafter requires some slight elaborations which we now describe.  To introduce path objects in the unpointed category $\cC$, we first choose a Reedy trivial fibration $\check{e}^\bullet\to e$ with $\check{e}^\bullet$ cofibrant in $\cC$.  Then, for an object $F^\bullet\in \cC$ with a map $\alpha\co \check{e}^\bullet\to F^\bullet$, we let $P_\alpha F^\bullet\in \cC$ be the \emph{path object} given by 
$$
P_\alpha F^\bullet~=~\hom^c(\Delta^1,F^\bullet)\times_{F^\bullet}\check{e}^\bullet~=~
\hom^c(\Delta^1,F^\bullet)\times_{F^\bullet\times F^\bullet}(\check{e}^\bullet\times F^\bullet)
$$
with the natural maps  $\check{e}^\bullet\to P_\alpha F^\bullet\to F^\bullet$  factoring $\alpha$.  We now replace Lemma \ref{lem:3.24} by the following: 

\begin{lem}
\label{lem:12.5}
For a $\G$--fibrant object $F^\bullet\in\cC$ with a map $\alpha\co \check{e}\to F^\bullet$, the natural map $P_\alpha F^\bullet\to F^\bullet$ (resp.\ $P_\alpha F^\bullet\to e$) has the right lifting property for $\G$--trivial cofibrations (resp.\ $\G$--cofibrations) in $\cC$.
\end{lem}

\begin{proof}
This follows easily from Lemma \ref{lem:3.23} since the map $\check{e}\to e$ has the right lifting property for $\G$--cofibrations.
\end{proof}

We likewise replace Lemma \ref{lem:3.25} by the following:

\begin{lem}
\label{lem:12.6}
If $F^\bullet\to e$ is a $\G$--trivial fibration with a $\G$--trivial cofibration $\alpha\co \check{e}^\bullet\to F^\bullet$, then $F^\bullet\to e$ has the right lifting property for $\G$--cofibrations.
\end{lem}

\begin{proof}
The $\G$--fibration $P_\alpha F^\bullet\to F^\bullet$ has a cross-section since it has the right lifting property for the $\G$--trivial cofibration $\alpha\co \check{e}^\bullet\to F^\bullet$ by Proposition \ref{prop:3.17}.  Hence $F^\bullet\to e$ has the right lifting property for $\G$--cofibrations since $P_\alpha F^\bullet\to e$ does by Lemma \ref{lem:12.5}.
\end{proof}

We now retain Lemma \ref{lem:3.26} but replace Proposition \ref{prop:3.27} by the following:

\begin{prop}
\label{prop:12.7}
A $\G$--trivial fibration $f\co X^\bullet\to Y^\bullet$ in $\cC$ has the right lifting property for $\G$--cofibrations.
\end{prop}

\begin{proof}
First suppose $X^\bullet$ is cofibrant.  By Proposition \ref{prop:3.21}, the map 
$X^\bullet\coprod\check{e}^\bullet\to e$ factors into a $\G$--cofibration   $\phi\co X^\bullet\coprod\check{e}^\bullet\to F^\bullet$ and a $\G$--trivial fibration $F^\bullet\to e$, and the map  
$(f,\phi)\co X^\bullet\to Y^\bullet\times F^\bullet$ factors into a Reedy cofibration $X^\bullet\to E^\bullet$ and a Reedy trivial fibration $E^\bullet\to Y^\bullet\times F^\bullet$.  Then the map $E^\bullet\to Y^\bullet$ is a $\G$--trivial fibration with the right lifting property for $\G$--cofibrations by Lemmas \ref{lem:3.22} and \ref{lem:12.6}.  The proof now proceeds as in \ref{prop:3.27}.
\end{proof}

We retain Proposition \ref{prop:3.28}, and thereby complete the proof of Theorem \ref{thm:12.4}.

\subsection{The unpointed theory}
\label{sec:12.8}
Our main definitions and results pertaining to $\G$--resolutions, right derived functors, and $\G$--completions in Sections \ref{sec:4}--\ref{sec:9} are now easily generalized to an unpointed model categories.  However, the main results in Sections \ref{sec:10}--\ref{sec:11} must be slightly modified since the $\G$--flatness condition for smash products (Definition \ref{defn:10.3}) must be replaced by a suitable $\G$--flatness condition for ordinary products.  This is easily accomplished when $\C=\Scal$ and, more generally, when the functor $-\times B\co \Ho\C\to\Ho\C$ has a right adjoint $(-)^B\co \Ho\C\to\Ho\C$ with $(\Omega^n Y)^B\homeq\Omega^n(Y^B)$ for each $B\in\Ho\C$ and $Y\in(\Ho\C)_*$.  We finally consider a general example leading to unpointed Bendersky--Thompson completions.  

\subsection{A general unpointed example}
\label{sec:12.9}
Let $\C$ be a left proper model category with a class $\Hcal$ of injective models in the associated pointed homotopy category $\Ho\C_*$.  As in \ref{sec:4.8}, the forgetful functor $J\co \Ho\C_*\to\Ho\C$ now carries $\Hcal$ to a class $J\Hcal$ of injective models in $\Ho\C$, and we obtain simplicial model categories $\cC^{J\Hcal}$ and $\cC^{\Hcal}_*$ together with Quillen adjoints $\cC^{J\Hcal}\leftrightarrows\cC^{\Hcal}_*$.  For an object $A\in\C_*$ with an $\Hcal$--resolution $A\to\bar{A}^\bullet$ in $\cC_*$, we easily deduce that $A\to\bar{A}^\bullet$ represents a weak $J\Hcal$--resolution of $A$ in $\cC$.  Thus, when $\C$ is bicomplete and simplicial, the $\Hcal$--completion $\hat{L}_{\Hcal}A\in\Ho\C_*$ represents the $J\Hcal$--completion $\hat{L}_{J\Hcal}A\in\Ho\C$, and we may view $\hat{L}_{J\Hcal}$ as an unpointed version of $\hat{L}_{\Hcal}$.

\subsection{The unpointed Bendersky--Thompson completions}
\label{sec:12.10}
The above discussion applies to give unpointed versions of the Bendersky--Thompson $E$--completions for ring spectra $E$ (\ref{sec:7.8}) and of the $p$--adic $K$--completion (\ref{sec:11.2}).

\end{document}